\documentclass[a4paper, reqno]{amsart}

\usepackage{amscd,amsthm,amsfonts,amssymb,amsmath,esint}
\usepackage{amsmath,tikz-cd}
\usepackage[colorlinks=true]{hyperref}
\usepackage{fullpage}
\usepackage[all]{xy}
\usepackage{mathtools}
\usepackage{mathrsfs}

    \newcommand{\BA}{{\mathbb {A}}} 
    \newcommand{\BC}{{\mathbb {C}}} 
     \newcommand{\BF}{{\mathbb {F}}}
    \newcommand{\BG}{{\mathbb {G}}}

     \newcommand{\BN}{{\mathbb {N}}}
     
    \newcommand{\BQ}{{\mathbb {Q}}} \newcommand{\BR}{{\mathbb {R}}}

     \newcommand{\BZ}{{\mathbb {Z}}}

    \newcommand{\CA}{{\mathcal {A}}} 
     
    \newcommand{\CE}{{\mathcal {E}}} 
     \newcommand{\CH}{{\mathcal {H}}}
     
     \newcommand{\CL}{{\mathcal {L}}}
     
    \newcommand{\CO}{{\mathcal {O}}} 
    \newcommand{\CQ}{{\mathcal {Q}}} 
    \newcommand{\CS}{{\mathcal {S}}}

    \newcommand{\alg}{{\mathrm{alg}}}

    \newcommand{\Seq}{{\mathrm{Seq}}}
    
    \newcommand{\corank}{{\mathrm{corank}}}
    
    \newcommand{\Cl}{{\mathrm{Cl}}}

    \newcommand{\diag}{{\mathrm{diag}}}

    \newcommand{\Gal}{{\mathrm{Gal}}} \newcommand{\GL}{{\mathrm{GL}}}

    \renewcommand{\Im}{{\mathrm{Im}}}

     \newcommand{\JL}{{\mathrm{JL}}}
    \newcommand{\Ker}{{\mathrm{Ker}}}

    \newcommand{\ord}{{\mathrm{ord}}} \newcommand{\rank}{{\mathrm{rank}}}
    \newcommand{\PGL}{{\mathrm{PGL}}}

    \renewcommand{\mod}{\ \mathrm{mod}\ }\renewcommand{\Re}{{\mathrm{Re}}}

    \newcommand{\Res}{{\mathrm{Res}}}\newcommand{\OGr}{{\mathrm{OGr}}}
    \newcommand{\Sel}{{\mathrm{Sel}}}

    \newcommand{\SL}{{\mathrm{SL}}}
     
    \newcommand{\SO}{{\mathrm{SO}}}\newcommand{\Sp}{{\mathrm{Sp}}}

    \newcommand{\tr}{{\mathrm{tr}}}\newcommand{\tor}{{\mathrm{tor}}}

\newcommand{\matrixx}[4]{\begin{pmatrix}
#1 & #2 \\ #3 & #4
\end{pmatrix} }        

    \font\cyr=wncyr10

    \newcommand{\Sha}{\hbox{\cyr X}}\newcommand{\wt}{\widetilde}
    \newcommand{\wh}{\widehat}
    
    \newcommand{\pair}[1]{\langle {#1} \rangle}

    \newcommand{\lra}{\longrightarrow}
    \newcommand{\ra}{\rightarrow} 
    \newcommand{\bs}{\backslash}

    \theoremstyle{plain}
    \newtheorem{thm}{Theorem}[section] \newtheorem{cor}[thm]{Corollary}
      \newtheorem{prop}[thm]{Proposition}
    \newtheorem {conj}[thm]{Conjecture} \newtheorem{defn}[thm]{Definition}
    
 \newtheorem{def-prop}[thm]{Definition-Proposition}

\theoremstyle{remark} \newtheorem{remark}[thm]{Remark}
\theoremstyle{remark} 
\theoremstyle{remark} 

    \numberwithin{equation}{section}

    \newcommand{\Prob}{{\mathrm{Prob}}}
    
\begin{document}

\title{The even parity Goldfeld conjecture: congruent number elliptic curves}

\author{Ashay Burungale and Ye Tian}

\address{Ashay A. Burungale:  California Institute of Technology,
1200 E California Blvd, Pasadena CA 91125}
\email{ashayburungale@gmail.com}

\address{Ye Tian: Academy of Mathematics and Systems
Science, MCM, HLM,  Chinese Academy of
Sciences, Beijing 100190,  and School of Mathematical Sciences, University of Chinese Academy of Sciences, Beijing 10049}
\email{ytian@math.ac.cn}

\begin{abstract}
In 1979 Goldfeld conjectured: 50\% of the quadratic twists of an elliptic curve defined over the rationals have analytic rank zero.
In this expository article we present a few recent developments towards the conjecture, especially its first instance - the congruent number elliptic curves.

\end{abstract}
\maketitle
\tableofcontents

\section{Introduction}

Representation of integers by ternary quadratic forms has rich history, yet it continues to be alluring. Sometimes, it is closely related to the arithmetic of quadratic twist family of elliptic curves defined over the rationals.

A positive square-free integer is called a congruent number if it is the area of a right triangle with rational side lengths. An important open problem: to determine whether or not a given integer is a congruent number, perhaps one of the oldest open problems (cf. \cite{T}). It is closely related to studying rational points on a certain quadratic twist family of elliptic curves - the congruent number elliptic curves.

We begin with
the Birch and Swinnerton-Dyer (BSD) conjecture for the congruent number elliptic curves in the guise of:
\begin{conj}\label{CN} Let $n$ be a positive square-free integer. The following are equivalent.
\begin{enumerate}
	\item[(a)] $n$ is a congruent number.

\item[(b)] Let $a=1$ if $2\nmid n$ and $a=2$ otherwise.  Let $\Sigma(n)$ be the set of  integral solutions to the equation
\[ \displaystyle{2ax^2+y^2+8z^2=\frac{n}{a}}. \]
Then,
$\#\left\{(x,y,z)\in \Sigma(n):  2|z\right\}=\#\left\{(x,y,z)\in \Sigma(n):  2\nmid z\right\}.$

\end{enumerate}
	
\end{conj}

Define
\begin{equation}\label{TunQ}
\CL(n)=\#\left\{(x,y,z)\in \Sigma(n):  2|z\right\}-\#\left\{(x,y,z)\in \Sigma(n):  2\nmid z\right\}.
\end{equation}
The non-vanishing of $\CL(n)$ may be determined in a finite number of steps, while an algorithm to determine whether a given $n$ is a congruent number remains elusive.

In view of Tunnell's theorem and the Coates--Wiles theorem: if $\CL(n)\neq 0$, then $n$ is not a congruent number (cf. \cite{CoWi}, \cite{Tu}). Conjecture \ref{CN} predicts the converse. One may ask:  
\begin{equation}\label{Q}\tag{Q}
\text{How often is $\CL(n)\neq0$?}
\end{equation}

\subsection{Main result} Our recent result \cite{BuTi2}: 
\begin{thm}\label{main}
  For a density one subset of positive square-free integers $n\equiv 1,2,3\pmod 8$,
  \[\CL(n)\neq 0.\]
\end{thm}
\begin{remark}
\noindent
\begin{itemize}
\item[$\circ$] A priori, an independent assertion: for a density one subset of $n\equiv 1,2,3\pmod 8$, $n$ is not a congruent number (cf. \cite{Sm1}).
\item[$\circ$] For $n\equiv 5,6,7\pmod{8}$, notice $\CL(n)=0$.
Conjecture \ref{CN} predicts that these $n$ are congruent. Over the last decade, arithmetic of Heegner point as pioneered by Heegner \cite{Hee}, \cite{Mon90} has led to a progress: \cite{Tian}, \cite{TYZ}, \cite{Sm1}.  It is now known that more than $50\%$ square-free positive integers $n\equiv 5,6,7\pmod{8}$ are congruent numbers (cf. \cite{TYZ}, \cite{Sm1}).

\end{itemize}
\end{remark}

\subsubsection{Congruent number elliptic curves}
Theorem \ref{main} yields the first instance of the influential (even parity) Goldfeld conjecture \cite{Go},  which concerns the distribution of analytic ranks in the quadratic twist family of elliptic curves over the rationals:

The congruent number problem may be rephrased in terms of the arithmetic of quadratic twist family of the congruent number elliptic curves\footnote{Note that $n$ is a congruent number if and only if $\rank_{\BZ}E^{(n)}(\BQ)>0$.}
\[E^{(n)}:ny^2=x^3-x.\]
Let $L(s,E^{(n)})$ denote the Hasse--Weil $L$-function of $E^{(n)}$. The integer $\CL(n)$ is closely related to the special $L$-value $L(1, E^{(n)})$.



\subsubsection{Outline}\label{ss:ot}


Theorem \ref{main} is a consequence of the following.
\begin{itemize}
\item An explicit Shimura--Shintani--Waldspurger correspondence \cite{Tu}: $$\CL(n)\neq 0\quad  \Longleftrightarrow\quad  L(1,E^{(n)})\neq 0.$$
  \item A $p$-converse theorem \cite{BuTi2}:
  \begin{equation}\label{$p$-cv}\tag{$p$-cv}
  \text{For any prime $p$, $\# \Sel_{p^\infty}(E^{(n)}/\BQ)<\infty \implies L(1,E^{(n)})\neq 0$,}
  \end{equation}
  where $\Sel_{p^\infty}(E^{(n)}/\BQ)$ denotes the $p^\infty$-Selmer group.
  \item A key progress towards Selmer-counterpart of the Goldfeld conjecture \cite{Sm1}:
  \[\Prob\left(\#\Sel_{2^\infty}(E^{(n)}/\BQ)<\infty\ \Big|\ n\equiv 1,2,3\pmod{8}\ \text{positive square-free}\right)=100 \%.\]
\end{itemize}

Our essential contribution is the $p$-converse theorem, especially for the prime $p=2$.  Now, an equivalent form of Theorem \ref{main}: the even parity Goldfeld conjecture for the congruent elliptic curves -
$$
\text{For a density one subset of positive square-free integers $n\equiv 1,2, 3 \mod{8}$, one has $L(1,E^{(n)})\neq 0$.}
$$

\begin{remark}
Since its proposal, the Goldfeld conjecture has been studied via diverse tools, yet an example remained elusive. Perhaps enigmatically the first example turns out to be the classical congruent number family. Time and again, the congruent number curves have influenced
the arithmetic of general elliptic curves over $\BQ$. Even a key precursor to \cite{Sm1} - the congruent number family \cite{TYZ}, \cite{Sm2}.
\end{remark}
\subsection{Plan} The article is essentially an elaboration of \S\ref{ss:ot}. It also reports on a generalisation of Tunnell's theorem to general quadratic twist families of elliptic curves \cite{HTX} and a preliminary investigation of a missing case in Smith's work \cite{FLPT}. The article is not meant as a survey. For instance, in view of \cite{BSTsv}, even the discussion of \eqref{$p$-cv} is succinct.

The text begins with the Goldfeld conjecture in \S\ref{s:Gd}.
Then \S\ref{s:Tug} presents a recent interrelation among ternary quadratic forms and central $L$-values of a quadratic twist family of elliptic curves over the rationals - a generalisation of Tunnell's theorem (the case of congruent number elliptic curves). Next an update of $p$-converse theorems appears in \S\ref{s:pcv}. Then \S\ref{s:dsel} briefly recalls a few conjectures regarding the distribution of Selmer groups associated to elliptic curves over a fixed number field and Smith's main result. Finally, \S\ref{s:dex} presents an exploratory study of a missing case in \cite{Sm2}, \cite{Sm1}.
\subsubsection*{Acknowledgement}
It is a pleasure to thank John Coates, Wei He, Shinichi Kobayashi, Jinzhao Pan, Dinakar Ramakrishnan, Alex Smith, Richard Taylor and Wei Zhang for helpful discussions and instructive comments. The authors cordially thank Chris Skinner and Shou-Wu Zhang for inspiring conversations. The article owes its existence to a generous suggestion of Dorian Goldfeld. The authors would like to express their sincere gratitude to Dorian Goldfeld also for his enticing conjecture.

A. B. is partially supported by the NSF grant DMS \#2001409, and Y. T. by the NSFC grants \#11688101 and \#11531008.

\section{Goldfeld's conjecture}\label{s:Gd}

\subsection{Backdrop}
\subsubsection{The set-up} An elliptic curve over the rationals is given by a projective curve with affine equation:
$$
A: y^{2}=x^{3}+ax+b
$$
for $a$, $b \in \BZ$ with $\Delta:=4a^{3}+27b^{2}\neq 0$.

The associated Hasse--Weil $L$-function
$L(s,A)$ is defined as an Euler product
 \[L(s,A):=\prod_{ \text{$p$ a prime}}L_{p}(p^{-s})^{-1}\]
 for $s\in\BC$, where
 $$
 L_{p}(X)=1-a_{p}X+pX^{2}, \quad a_{p}=p+1-\#A(\BF_p)
 $$
 for $p\nmid 2\Delta$. Define
 $$
 \Lambda(s,A):=N^{s/2}\cdot 2(2\pi)^{-s}\Gamma(s)L(s,A)
 $$
 for $N$ the conductor.

 In view of the Hasse bound $|a_{p}|\leq 2\sqrt{p}$,
the Euler product is absolutely convergent for $\Re(s)>3/2$. The elemental modularity:

\begin{thm}
The Hasse--Weil $L$-function $L(s,A)$ has entire continuation,
which satisfies the functional equation
$$
\Lambda(s,A)=\varepsilon(A)\Lambda(2-s,A),
$$
where $\varepsilon(A)\in\{\pm 1\}$ denotes the root number.
\end{thm}
The central vanishing order - $\ord_{s=1}L(s,A)$ - is referred to as the analytic rank of $A$.

\subsubsection{The Birch and Swinnerton-Dyer conjecture}
\begin{conj}[The BSD conjecture]\label{cBSD}
Let $A$ be an elliptic curve over $\BQ$.
\begin{itemize}
  \item[(a)] $\ord_{s=1}L(s,A)=\rank_{\BZ}A(\BQ)$
  \item[(b)] The Tate--Shafarevich group $\Sha(A/\BQ)$ is finite and
   $$\frac{L^{(r)}(1,A)}{r!\cdot \Omega_A \cdot R_A}=\frac{\prod_{\ell}c_\ell(A)\cdot  \#\Sha(A/\BQ)}{\# A(\BQ)_{\tor}^2}$$
  for $r=\ord_{s=1}L(s,A)$,  $\Omega_{A} \in \BC^{\times}$ the N\'eron period, $R_A$ the regulator of the N\'eron--Tate height pairing on $A(\BQ)$ and $c_{\ell}(A)$ the Tamagawa number at $\ell$.
  \end{itemize}

\end{conj}

The Tate--Shafarevich group $\Sha(A/\BQ)$ is defined as
 \[\Sha(A/\BQ)=\Ker\left(H^1(\BQ,A)\ra \prod_p H^1(\BQ_p,A)\right).\]
It may be interpreted as the isomorphism classes of $A$-torsors $C$ such that $C(\BQ_p)$ is non-empty for all primes $p$.


\begin{remark}
For a brief introduction, one may refer to the recent survey \cite{BSTsv}.
\end{remark}
\subsection{Goldfeld's conjecture}
An individual invariant may often be delicate to study, an emerging theme is to instead investigate its  variation in a family\footnote{which may shed some light on the individual invariant}. In the late 1970's Goldfeld pioneered the exploration of quadratic twist families of elliptic curves over the rationals.

Let $A:y^{2}=x^3+ax+b$ be an elliptic curve over the rationals as above. For a square-free integer $d$, consider the quadratic twist $A^{(d)}: dy^{2}=x^3+ax+b.$
A principal insight of Goldfeld is that the underlying analytic or arithmetic invariants often vary systematically in the quadratic twist family $\{A^{(d)}\}_{d}$.

\subsubsection{The conjecture}
In 1979 Goldfeld \cite{Go} proposed the following
\begin{conj}[Goldfeld's conjecture] \label{Gd}Let $A$ be an elliptic curve over $\BQ$.

Then, for a density one subset of square-free integers $d$ with $\varepsilon(A^{(d)})=+1$ (resp.~$\varepsilon(A^{(d)})=-1):$
\[
\ord_{s=1}L(s,A^{(d)})=0
,\quad \text{(resp.~$\ord_{s=1}L(s,A^{(d)})=1
$)}.\]

\end{conj}

We refer to the sign $+1$ (resp.~$-1$) part as the even (resp.~odd) parity Goldfeld conjecture. It may be easily seen that $50\%$ of the quadratic twists have sign $\pm 1$.
\begin{remark}
The core of Goldfeld's conjecture is his minimalist principle: Often for natural families of elliptic curves over
$\BQ$ - not just the quadratic twist families - the subfamily with root number $+1$ (resp.~$-1$) has generic analytic rank $0$ (resp.~$1$). In particular, the distribution of analytic rank is the same as that of the root number.
\end{remark}
\begin{remark}\label{Gdg}
\noindent
It is natural to seek an analogue of the conjecture over number fields.
In general, the root number variation in a quadratic twist family may notably differ.
\begin{itemize}
\item[$\circ$] A counterpart over number fields: \cite[Conj. 7.12]{KMR}.
\item[$\circ$] One may also seek a variant of the conjecture for a self-contragredient cuspidal automorphic representation of $\GL_{2}(\BA_F)$ for $F$ a number field.
Such an investigation appears in \cite{Ba}. Also see Conjecture \ref{Gdv} below.
\item[$\circ$] An instance of a contrasting root number variation: Let $E/F$ be an elliptic curve with everywhere good reduction, $F$ with no real places but odd (resp. even) number of complex places. Then the root number is given by $\varepsilon(E)=-1$ (resp. $\varepsilon(E)=+1$), further any quadratic twist of $E$ also has root number $-1$ (resp. $+1$). Such examples perhaps first appeared in \cite{DD}. Over $\BQ(\sqrt[6]{-11})$, the elliptic curve
 \[y^2=x^3+\frac{5}{4}x^2-2x+7\] has everywhere good reduction and
  its any quadratic twist $E'$ satisfies $\varepsilon(E')=-1$. In contrast, over $\BQ(\sqrt[4]{-37})$, the elliptic curve \[y^2=x^3+x^2-12x-\frac{67}{4}\]has everywhere good reduction and its any quadratic twist $E'$ satisfies $\varepsilon(E')=+1$.
\end{itemize}
\end{remark}
\subsubsection{An existence} We recall a mild, yet general result towards Conjecture \ref{Gd} (cf. \cite{BFH}).

Let $F$ be a number field and $\pi$ a self-contragredient cuspidal automorphic representation of $\GL_{2}(\BA_F)$. Then its root number $\varepsilon(\pi)\in \{\pm 1\}$ satisfies
\[(-1)^{\ord_{s=1/2}L(s,\pi)}=\varepsilon(\pi),\]
where $L(s,\pi)$ is the $L$-function of $\pi$. Any quadratic twist of $\pi$ is also self-contragredient.
\begin{thm}\label{exs}
  Let $\varepsilon\in\{\pm1\}$ and $\chi$ be a quadratic character over $F$ such that $\varepsilon(\pi\otimes\chi)=\varepsilon$. Let $S$ be a finite set of places of $F$.

  Then, among the quadratic characters $\chi'$ over $F$ with $\chi'_v=\chi_v$ for $v\in S:$ there exist infinitely many $\chi'$ such that $\ord_{s=1/2} L(s, \pi\otimes \chi')=0$ (resp.~$1$) if $\varepsilon=+1$ (resp.~$-1$).

\end{thm}
\subsubsection{The conjecture, again}
In light of Goldfeld's minimalist principle, one may naturally propose:

\begin{conj}\label{Gdv}
 Let $\pi$ be a self-contragredient cuspidal automorphic representation of $\GL_{2}(\BA_F)$.  Let $\chi$ be a quadratic character over $F$ with $\varepsilon(\pi\otimes \chi)=\varepsilon$ and let $S$ be a finite set of places of $F$.

Then, among the quadratic characters $\chi'$ over $F$ with
\begin{itemize}
\item[(i)] $\varepsilon(\pi\otimes \chi')=\varepsilon$  and
\item[(ii)] $\chi'_v=\chi_v$ for $v\in S$,
\end{itemize} the density of $\chi'$ with $\ord_{s=1/2} L(s, \pi\otimes \chi')=0$ (resp. $1$) is one if $\varepsilon=+1$ (resp. $-1$).
\end{conj}

\section{Tunnel's theorem, generalised}\label{s:Tug}
The section reports on a recent  generalisation \cite{HTX} of Tunnell's theorem to general quadratic twist families of elliptic curves over $\BQ$ (cf. Theorem \ref{thmm}) . The strategy - a departure from Tunnell's method - employs  general explicit Waldspurger formula \cite{CST} and  explicit theta liftings.

\subsubsection{Notation}  For $n\in \BQ^\times$, let $\chi_n$ be the quadratic character over $\BQ$ corresponding to the extension $\BQ(\sqrt{n})$.  For $N\in \BZ$ a positive integer, $\chi: (\BZ/N\BZ)^\times \ra \BC^\times$ a Dirichlet character,  and  $k\in\frac{1}{2}\BZ$ such that  $4|N$ if $k\notin\BZ$, let  $M_{k}(N,\chi)$ (resp.~$S_k(N,\chi)$) denote the space of modular forms (resp.~cusp forms)
of weight $k$, level $\Gamma_0(N)$, and character $\chi$. These spaces are endowed with Hecke action. 

\subsection{Tunnell's theorem}
\subsubsection{The theorem} Let $E^{(n)}: y^2=x^3-n^2x$ be the congruent elliptic curve, where $n$ is a positive square-free integer.

A link among the central $L$-values and ternary quadratic forms:
\begin{thm}[Tunnell's theorem]\label{Tun}
There are weight $3/2$ modular forms,
$$
\sum_{n=1}^\infty a_nq^n\in S_{3/2}(128,\mathbf{1}),\quad
\sum_{n=1}^\infty b_nq^n\in S_{3/2}(128,\chi_2)
$$
such that for all positive square-free integers $n$,

\[\CL(n)= \begin{cases}
a_n, \\
b_{n/2},
\end{cases},\frac{L(1,E^{(n)})}{\Omega/\sqrt n}= \CL(n)^2\cdot \begin{cases}
\frac{1}{16},&\quad\quad  \text{if $2\nmid n$} \\
\frac{1}{8},&\quad \quad \text{if $2\mid n$}
\end{cases}.
\]Here $\CL(n)$ as in \eqref{TunQ}
and $\Omega=\displaystyle{\int_1^\infty\frac{\mathrm dx}{\sqrt{x^3-x}}}$. (cf. \cite{Tu}, \cite{Qin})
\end{thm}

In light of Tunnell's theorem, the central $L$-values of the quadratic twist family of congruent elliptic curves are modular. Furthermore, the theorem gives an effective way to compute the $L$-values.

\subsubsection{Tunnell's proof}

 The key tool is a fundamental theorem of Waldspurger, which connects
\begin{itemize}
\item The Fourier coefficients of half weight modular forms that are Shimura equivalent to a given elliptic newform $\varphi$,
\item The central $L$-values of the quadratic twists of $\varphi $.
    \end{itemize}
\vskip2mm
{\it{Shimura equivalence}.}
The Shimura equivalence connects - weight $2$ and weight $3/2$ - Hecke eigenforms.

Given a newform $\varphi \in S_2(M,\chi^2)$ and an positive integer $N\in 4\BZ\cap 2M\BZ$, the subspace of $S_{3/2}(N,\chi)$ Shimura equivalent to $\varphi $ is given by
\[S_{3/2}(N, \chi, \varphi):=\left\{f\in S_{3/2}^{\perp}(N,\chi)\ \big|\ T_{p^2}f=a_p(\varphi)f\ \text{for all $p\nmid N$}\right\}.\] Here $S_{3/2}^{\perp}(N,\chi)$ is the subspace of $S_{3/2}(N,\chi)$ orthogonal to one variable theta series.

Let  $\varphi \in S_{2}(M, \chi^2)$ be a newform and $\pi=\otimes_v \pi_v$  the irreducible automorphic representation of $\GL_2(\BA)$ associated to $\varphi $.  After Flicker, there exists an integer $N$ such that
$S_{3/2}(N, \chi, \varphi)\neq 0$ if and only if the following hypothesis holds:
 If $\pi_v=\pi(\xi_{1,v},\xi_{2,v})$ is a principal series with associated characters  $\xi_{1,v}$, $\xi_{2,v}$, then
 \begin{equation}\label{H}\tag{H}\xi_{1,v}(-1)=\xi_{2,v}(-1)=1.
 \end{equation}

\begin{thm}[Waldspurger]\label{cor} Let  $\varphi \in S_{2}(M, \chi^2)$ be a newform that satisfies the hypothesis \eqref{H}. Let $f=\sum a_n q^n\in S_{3/2}(N, \chi, \varphi)$ with $N\in 4\BZ \cap 2M \BZ$.

If $n_1$, $n_2$ are positive square-free integers with $n_1/n_2\in \BQ_p^{\times 2}$ for all $p|N$, then

\[a_{n_1}^2\cdot L\left(1, \varphi \otimes \chi_0^{-1}\chi_{n_2}\right)\sqrt{n_2}\chi(n_2/n_1)=a_{n_2}^2\cdot L\left(1, \varphi \otimes \chi_0^{-1}\chi_{n_1}\right)\sqrt{n_1}. \]
Here $\chi_0(n)=\chi(n)\left(\frac{-1}{n}\right)$. (cf. \cite{Wal})

\end{thm}
\begin{remark}A consequence: Assume that $a_n\neq 0$  for a positive square-free $n$.  Then $L(1, \varphi \otimes\chi_0^{-1}\chi_n)\neq 0$, moreover,  for any positive square-free $m$ with $m/n\in \BQ_p^{\times 2}$ for all $p|N$, one has
$$a_m\neq 0 \iff L(1, \varphi \otimes\chi_0^{-1}\chi_m)\neq 0.$$
\end{remark}

Let $\varphi $ be the newform associated to $E: y^2=x^3-x$ of weight $2$ and level $32$.
Now, consider inverse image of the Shimura equivalence for $\varphi $.
As $S_{3/2}(64,\chi)$ with $\chi^2=1$ is
generated by one variable theta series,
one may resort to $S_{3/2}(128,\chi)$, which is interlaced with weight $1/2$ modular forms.
Indeed, $S_{3/2}(128,\chi)$ and $M_{1/2}(128,\chi)$ are $3$-dimensional, the latter generated by one variable theta series, so multiplication by the unique newform in $S_{1}(128,\chi_{-2})$ induces an explicit isomorphism
 \[M_{1/2}(128,\chi\chi_2)\stackrel{\sim}{\lra}  S_{3/2}(128,\chi).\]
This gives rise to $f_1\in S_{3/2}(128,\mathbf{1})$ and $f_2\in S_{3/2}(128,\chi_2)$, which are Shimura equivalent to $\varphi $ satisfying
$$\text{$a_{1}(f_1)$, $a_3(f_1)$, $a_1(f_2)$ and $a_5(f_2)$ are non-zero. }$$

Then Theorem \ref{Tun} is just a special case of Theorem \ref{cor}.

\begin{remark}
If the genus class of a definite integral ternary quadratic form consists of only two forms, then the difference of the associated theta series is an eigen cusp form of weight $3/2$.
Qin found the above weight $3/2$ modular forms as such (cf. \cite{Qin}).

\end{remark}

\subsection{General counterpart}
One may seek a generalisation of Tunnell's theorem:
As $A$ varies in the quadratic twist family $\CE$ of an elliptic curve over $\BQ$ -
\begin{itemize}
 \item[(a)]  Encapsulate modularity of a certain ``square root'' of \[L^{\alg}(1,A):=\frac{L(1,A)}{\Omega_A}\in \BQ.\]
     \item[(b)] Offer an effective algorithm to compute $L^{\alg}(1,A)$ in terms of ternary quadratic forms. In particular, an algorithm to determine non-vanishing\footnote{In turn to determine positivity of $\rank_{\BZ}A(\BQ)$ in finitely many steps as
 $L(1,A)\neq 0\Rightarrow \rank_{\BZ}A(\BQ)=0$ (cf. \cite{CoWi}, \cite{K}, \cite{Ko}, \cite{SU}).
 Assuming (the rank part of) the BSD conjecture, this is an effective algorithm.
} of $L(1,A)$.
     \end{itemize}

The first is essentially addressed by Waldspurger \cite{Wal}. The second has been studied extensively, a notable progress due to Gross \cite{Gr2} - quadratic twists of elliptic  curves with prime conductor - via Waldspurger formula for toric periods. (See also an extension \cite{BSP} of Gross' work to the square-free conductor case under some local conditions.)
Definite ternary quadratic forms are elemental to the approach.

In a recent joint work
\cite{HTX} of the second author, a Waldspurger-style result is reproven in the context of automorphic forms via theta lifting and Waldspurger formula for toric periods (via the local test vector theory developed in \cite{CST}). The approach leads to a generalization of the aforementioned results - due to Tunnell and Gross - for any quadratic twist family. 
In light of the explicit Waldspurger formula \cite{CST},
(b) is now available
for the general quadratic twist family of weight $2$ newforms with trivial character.

In general, due to potential local obstruction arising from the action of Atkin--Lehner operators,
it is essential to consider a particular subset of the integral solutions of relevant ternary quadratic form. We call them oriented solutions (cf. \eqref{ort}).
The previous results toward $(b)$ implicitly assume the vanishing of the local obstruction (cf. Remark \ref{rem1}).

In addition, perhaps surprisingly, it is crucial\footnote{If $\CE$ is a family of CM curves or there exists $A\in\CE$ such that the conductor of $A$ is not a square, then definite ternary quadratic forms suffice.} to resort to certain indefinite ternary quadratic forms (cf. Remark \ref{rem1}).
\subsubsection{Theta lifting}
 The theory of theta lifting generalizes the classical construction of half weight modular forms from quadratic forms.

Let $B/\BQ$ be a definite quaternion algebra. Then $V:=B^{\tr=0}$ is a quadratic space with quadratic form $q$ given by minus of the reduced norm. Let $H=\SO(V)=PB^\times$ and $\BG=\wt{\SL_2(\BA)}$ be the metaplectic double cover of $\SL_2(\BA)$. Fix a non-trivial additive character $\psi$ of $\BQ\bs\BA$.
There is a Weil representation $w$ (associated to $\psi$) of $H(\BA)\times \BG$ on $\CS(V(\BA))$.
Let $\CA_0(H)$ (resp. $\CA_{0}(\BG)$) be the space of automorphic forms on $H(\BA)$ (resp. $\BG$).
 Theta lifting (associated to $\psi$) is a systematic mechanism to construct automorphic forms on $\BG$ from $H$ and Schwartz functions $\CS(V(\BA))$, via the Weil representation.
For each $\phi\in \CS(V(\BA))$, the theta kernel function
\[\theta_\phi: (h,g)\mapsto \sum_{x\in V}(w(h,g)\phi)(x)\] is an automorphic form on $H(\BA)\times \BG$ and gives a ($H(\BA)\times \BG$)-equivariant map
\[\theta:\CA_0(H)\times\CS(V(\BA))\ra \CA_0(\BG),\quad (f,\phi)\mapsto \left(\theta_f^\phi:g\mapsto \int_{H(\BQ)\bs H(\BA)}f(h)\theta_\phi(h,g)dh\right).\]

 Let $\pi\subset \CA_{0}(H)$ be an irreducible cuspidal representation. Recall the theta lifting of $\pi$ is defined to be \[\theta(\pi):=\{\theta_{f}^\phi\ \big|\ f\in \pi,\phi\in \CS(V(\BA))\},\] which is an irreducible cuspidal automorphic representation of $\BG$.
It is known that $\theta(\pi)\neq 0$ if and only if $L(\frac{1}{2},\pi)\neq 0$.

 Let $\CA_{0, 2}(H)$ denote the space of cuspidal automorphic forms on $H(\BA)$ with the $H(\BR)$-action being trivial and let $\pi\subset \CA_{0,2}(H)$ be an irreducible. Our goal is to explore precise relation between Fourier coefficients of $\theta_f^\phi$ and the central $L$-values of quadratic twists of $\pi$, by choosing explicit test vectors $f\in \pi$ and $\phi\in \CS(V(\BA))$.

First of all, the Fourier coefficients of theta liftings $\theta_f^\phi$ have a natural connection with ternary quadratic forms and toric periods.

 For $m\in \BQ$, $f\in \CA_0(H)$ and $\phi\in \CS(V(\BA))$, define the $m$-th Fourier coefficient of $\theta_f^\phi$ by \[W_m(\theta_f^\phi)=\int_{\BQ\bs\BA}\theta_f^\phi\begin{pmatrix}
                                                          1 & x \\
                                                          0 & 1
                                                        \end{pmatrix}\psi(-mx)dx.\]

    \label{cpt}
Pick $f\in \pi$,  then $f$ is constant at infinity and fixed by an open compact subgroup $U$ of $H(\BA_f)$.
 Pick $\phi\in \CS(V(\BA))$ such that $\phi|_{V_m(\BA)}\in \CS(V_{m}(\BA))^{U}$, where $V_{m}(\cdot )=\{x\in V(\cdot)\ \big|\ q(x)=m\}$ for $\cdot=\BQ,\BQ_p,\BA$.
  Assume that $\phi=\phi_{\text{fin}}\otimes \phi_\infty$ with finite part given by $\phi_{\text{fin}}=\mathbf{1}_{{Z}}$ for an open compact set ${Z}\subset V(\BA_f)$ and that
 $\phi_\infty$ factors through the reduced norm. Say $Z=\wh{L}$ for a lattice $L$ in $V$. 
 
 On the one hand, $W_m(\theta_f^\phi)$ has expression in terms of ternary quadratic forms:
  \[\label{FT}\tag{FT}W_m(\theta_f^\phi)\doteq\sum_{[h]\in B^\times\bs \wh{B}^\times/U}\frac{f(h)}{w_h} \cdot  \# (V_m(\BQ)\cap Z^{h})  \] with $w_h=\# H(\BQ)\cap hUh^{-1}$ and $Z^h= h{Z}h^{-1}$. Here `$\doteq$' denotes equality up to an explicit non-zero constant depending on $\phi_\infty$ and choice of measure.

  On the other hand, using Witt theorem, $W_m(\theta_f^\phi)$ can be expressed in terms of toric period:
 \label{fourtor}For  $x\in V_m(\BQ),$
  \[\label{FP}\tag{FP}W_m(\theta_f^\phi)=\int_{T_x(\BA)\bs H(\BA)}\phi(h^{-1}\cdot x)P_{T_x}(f^h) dh.\]
Here $\cdot$ denotes the conjugate action of $H$ on $V$, $T_x\subset H$  is the stabilizer of $x$,  and $P_{T_x}(f)$ is the toric period
\[P_{T_x}(f):=\int_{T_x(\BQ)\bs T_x(\BA)}f(t)dt.\]

  For $f\in \pi$, $\phi\in \CS(V(\BA))$, we seek a relation between the Fourier coefficients $W_m(\theta_f^\phi)$ of $\theta_f^\phi$ and the quadratic twist central $L$-values $L(\frac{1}{2},\pi\otimes \chi_{m})$ as $m$ varies.
 Let $M$ be the conductor of the Jacquet--Langlands correspondence $\sigma$ of $\pi$. 
\begin{defn}
\label{M-equivalent} Let $C\subset \BQ^\times$ be an equivalence class with respect to the relation:
\[a\sim b \quad \iff \quad a/b\in\BQ_v^{\times 2}\quad\text{for all}\quad v\mid 2M\infty.\]
We call $C$ an $M$-equivalence class.
\end{defn}

The root number of $\sigma\otimes \chi_m$ is independent of $m\in C$, denoted by $\varepsilon(C)$. Assume $\varepsilon(C)=1$ and that $(\pi,C)$ satisfies the Tunnell--Saito condition (cf. \eqref{TS}). In the explicit Waldspurger formula \cite{CST} (see Theorem \ref{CST1}), as $m$ varies in $C$, one can find a common vector $f\in \pi$, such that the toric period $P_{T_x}(f)$,  is proportional to the  base change central $L$-value $L(\frac{1}{2},\pi_{K_m})=L(\frac{1}{2}, \pi)L(\frac{1}{2}, \pi\otimes \chi_m)$. Here $x\in V$ with $q(x)=m$ and $K_m=\BQ(\sqrt{m})$.

The following notion of
an admissible pair $(f,\phi)$ for $f\in\pi$ and $\phi\in \CS(V(\BA))$ is elemental in the context of Shimura equivalence.
 It is purely local, perhaps easier to check and more flexible than Shimura equivalence.

 From now on, further choose $\psi$ in the Weil representation with $\psi_\infty(x)=e^{2\pi i c x}$ for $c\in \BQ_{<0}$.

\begin{defn}\label{dt}
  Let $f=\otimes f_v\in \pi$  and $\phi=\otimes_v \phi_v\in \CS(V(\BA))$ be pure tensors. We call $(f,\phi)$ an admissible pair if
  \begin{itemize}
   \item for $p\nmid 2M\infty$,  $f_p$ is spherical and $\phi_p=\mathbf{1}_{M_2(\BZ_p)^{\tr=0}}$, 
  \item  $\phi_\infty(x)=e^{2\pi|c|q(x)}$, 
  \item for $p|2M$, $\phi_p\in \CS(V(\BQ_p))$ is invariant under left multiplication by $1+2p\BZ_p$.  \end{itemize}

\end{defn}

Define normalized Fourier coefficients of $\theta_f^\phi$: \[a_{n}(\theta_f^\phi):=e^{2\pi |c|n }W_{-n}(\theta_f^\phi),\quad n\in \BQ.\]
Relation of  explicit Waldspurger formulae with different test vector is given in \cite{CST}.  We  compare  $a_m(\theta_f^\phi)$ with the toric periods defined by another test vector in \cite{CST}.  It turns out, for an admissible pair $(f, \phi)$,  the ratio is constant as $m$ varies in an equivalence class.
 We have (cf. \cite{HTX}):
\begin{thm}\label{walds}
Let $\pi\subset \CA_{0,2}(H)$ with $L(\frac{1}{2},\pi)\neq 0$.
Let $(f,\phi)$ be an admissible pair.

 Then for any positive square-free integers $n_{1}, n_{2}$ with $n_1/n_2\in \BQ_p^{\times2}$ for all $p|2M$,
\[a_{n_1}^2(\theta_f^\phi)L\left(\frac{1}{2},\pi\otimes\chi_{-n_2}\right)\sqrt{n_2}=a_{n_2}^2(\theta_f^\phi)L\left(\frac{1}{2},\pi\otimes\chi_{-n_1}\right)\sqrt{n_1}.\]
\end{thm}
However, the above theorem is ineffective since it does not offer construction of an admissible pair such that the Fourier coefficient with a given index is non-zero.  Actually, the construction is the crux of the main result of \cite{HTX}, which relates quadratic twist $L$-values to ternary quadratic forms.

 \subsubsection{Quadratic twist subfamilies}
 Let $\CE$ be the quadratic twist family of an irreducible cuspidal automorphic representation of $\PGL_2(\BA)$ whose infinite component is the discrete series of weight $2$.  We first partition the family into finitely many subfamilies and then for each subfamily, we construct an admissible pair that is effective for this subfamily:
  the Fourier coefficients of its theta lifting
 interpret the central $L$-values of the representations in this subfamily. Furthermore, we also seek to relate the Fourier coefficients to ternary quadratic forms.

Let $\sigma\in \CE$ be such that $L(\frac{1}{2},\sigma)\neq 0$ and $M$ denote the conductor of $\sigma$. Let $C\subset \BQ^\times$ be an $M$-equivalence class. Then we have a subfamily of $\CE$ given by
\[\CE_{C}:=\{\sigma\otimes\chi_{m}\ \big|\ m\in C\}\subset \CE.\]  Actually,
$\CE$ can be covered by finitely many $\CE_{C}$'s ($\sigma$ may differ).

 Notice it suffices to consider a set-up: $L(\frac{1}{2},\sigma)\neq 0$ and $C$ an $M$-equivalence class with
 $\varepsilon(C)=1$. Let $B$ be the quaternion algebra over $\BQ$ such that $(\sigma,C)$ satisfies Tunnell--Saito condition, i.e. for each $K=\BQ(\sqrt{m})$, $m\in C:$
\[\varepsilon(\sigma_v,\mathbf{1}_{K_{v}})=\eta_v(-1)\epsilon(B_v) \label{TS}\tag{TS}\] holds for all places $v$ of $\BQ$. Here $\varepsilon(\sigma_v,\mathbf{1}_{K_v})$ is the local root number at $v$ of the Rankin--Selberg $L$-function $L(s,\sigma\times\mathbf{1}_{K})$, $\eta_v$ the quadratic character associated to the extension $K_v/\BQ_v$ and $\epsilon(B_v)\in \{\pm 1\}$ the invariant of $B$ at $v$.
Let $H=PB^\times$ and  $\pi\subset \CA_{0}(H)$ be with $\pi^{\JL}=\sigma$.

  For simplicity, we often assume $\sigma$ corresponds to an elliptic curve over $\BQ$ and that $B$ is definite\footnote{ In fact, to construct half weight modular forms effectively for $\CE_C$ and to connect with ternary quadratic forms, the ensuing approach  works well for a general quaternion algebra.}.


\subsubsection{Quadratic twist $L$-values and toric periods} We first present the explicit Waldspurger formula in \cite{CST} which connects toric periods to quadratic twist $L$-values.

  Let $\pi\subset \CA_{0,2}(H)$ be irreducible with conductor $M$. For simplicity, assume that $\pi$ corresponds to an elliptic curve over $\BQ$ via Jacquet-Langlands correspondence. Let $C\subset \BQ^\times$ be an $M$-equivalence class. Assume $(\pi, C)$ satisfies the Tunnell--Saito condition \eqref{TS}.

 Let $K\subset B$ with $K\simeq \BQ(\sqrt{m})$, $m\in C$. An order $R\subset B$ is called admissible (with respect to $(\pi,K)$) if its discriminant equals $M$ and $R\cap K=\CO_{K}$. Then the following space is one dimensional:
\label{CSTt}
\[V(\pi,K):=\left\{f\in \pi \ \Biggr|\ \substack{\displaystyle f  \ \text{is invariant under $\wh{R}^\times$ and ${K_{p}^\times}$ for any $p|M$ ramified in $K$}}
\right\}.\]We call its generator an admissible test vector for $(\pi,\mathbf{1}_K)$ with level $\wh{R}^\times$.
(A proof of the fact that the space $V(\pi,K)$ is one dimensional and consists of pure tensors appears in \cite{CST}, which generalizes the newform theory.)
\label{distl}
Let $\Sigma_C\subset C$ denote the subset of fundamental discriminants in $C$. A vector $f_0\in \pi$ is admissible for $C$ if there exists $x\in V_{D_0}(\BQ)$ for some $D_0\in \Sigma_C$ and an admissible order $R_0$ (with respect to $(\pi,\BQ(x))$) such that $f_0$ is an admissible test vector for $(\pi,\mathbf{1}_{\BQ(x)})$ with level $\wh{R}_0^\times$. Let $f_0$ be admissible for $C$, then for any $D\in \Sigma_C$ and $x_{D}\in V_{D}(\BQ)$, there exists $h_D\in H(\BA)$ such that $f_D:=f_0^{h_D}$ is an admissible test vector for $(\pi,\mathbf{1}_{\BQ(x_D)})$ with level $\wh{R}_D^\times$ for $R_D:=h_D\wh{R}_0 h_{D}^{-1}\cap B$ (an admissible order).

In particular, we may choose $h_{D_0}=1$ and then $f_0=f_{D_0}$. From now,   fix such choices.

For $D\in \Sigma_C$, let $K_D=\BQ(x_D)$ and $\CO_D$ the ring of integers of $K_D$. Denote \[P_{K_D}^0(f_D):=\sum_{t\in \Cl(K_D)}f_D(t).\]
\begin{thm}[Explicit Waldspurger formula] \label{CST1}Let $f_0$ be a non-zero $C$-admissible test vector, valued in the rationals. 

Then
  there exists a constant $k_C\in \BQ^\times$ dependent only on $\pi$ and $C:$ 
  \[\frac{\sqrt{D}\cdot L\left(\frac{1}{2},\pi_{K}\right)}{\Omega_\sigma^+\Omega_\sigma^-}=k_C\cdot \frac{|P_{K}^0(f_D)|^2}{[\CO_D^\times:\BZ^\times]^2}\in \BQ\]
  for all $D\in \Sigma_C$, $K=\BQ(\sqrt{D})$. Here $(\Omega_\sigma^+, \Omega_\sigma^-)\in (\BR^\times\times i\BR^\times)$ are Shimura's fundamental periods associated to $\sigma$. (cf.  \cite{CST})
  \end{thm}

\subsubsection{Optimal choices}

For $B$ and $f_0=f_{D_0}\in \pi$  a $\BZ$-primitive $C$-admissible test vector as in Theorem \ref{CST1}, one may seek a $\phi_0$ inherent to the ternary quadratic form such that exactly one toric period appears in the Fourier coefficients\footnote{Note that the $D$-th Fourier coefficient $W_D(\theta_f^\phi)$ only depends on the restriction $\phi|_{V_D(\BA)}$} of $\theta_{f_0}^{\phi_0}$.  One may even consider all $D\in \Sigma_C$ simultaneously.

To begin, choose $h_{D}$ for $D\in \Sigma_{C}$: for  $v| 2M\infty$ and $D,D'\in \Sigma_C$, \[h_{D',v}^{-1}\cdot(\sqrt{D/D'} x_{D'})=h_{D,v}^{-1}\cdot x_D\in V(\BQ_v).\] Here for $p|2M$, view $\sqrt{D/D'}$ as an element in $ (1+2p\BZ_p)$. A simple yet key result \cite{HTX}:
\begin{prop} There exists a lattice $L_0\subset V$ and an open compact subgroup $U_0\subset \wh{R}_0^\times$ such that  $$L_D:=h_D\wh{L}_0h_D^{-1}\cap V$$ is an $x_D$-distinguished lattice with level $U_{D}:=h_{D}U_0h_{D}^{-1}\subset \wh{R}_D^\times$ for any $D\in \Sigma_C$ - the definition being: \[L_{D, p}\cap V_{D}(\BQ_p)=U_{D,p}\cdot \{\pm x_D\}\]
for all prime $p$.
\end{prop}


For unique toric integral to appear in the Fourier coefficients,  it is essential to  isolate a certain symmetry arising from the abelian group structure of lattices, thereby switching to a finer structure: For a prime $p|2M$, define a local condition $L_{D,p}^o\subset L_{D,p}$ by
\[L_{D,p}^o=\Big{\{}\ell\in L_{D,p}\ \Big|\ \ell\in U_{D,p}\cdot((1+2p\BZ_p)x_{D}) \Big{\}}\]

Now let $\phi_D=\phi_{D,\text{fin}}\otimes e^{2\pi|c|q(\cdot)}$:  $\phi_{D,\text{fin}}$ the characteristic function of $\wh{L}_D^{(2M)}\cdot \prod_{p|2M}L_{D, p}^o$.

\begin{defn}
  The special admissible pairs $(f_D,\phi_D)$, $D\in \Sigma_C$ are called distinguished.
\end{defn}

\label{distp} The theta lifting \[\theta_C:=\theta_{f_D}^{\phi_D}\] does not depend on the choice of $D$.
Now let $(f,\phi,L,U,\{L_p^o\}_{p|2M})$ be any one of
$$\big(f_D, \phi_D, L_D, U_D, \{L_{D,p}^o\}_{p|2M}\big), \qquad D\in \Sigma_C. $$
\begin{defn}
For each $[h]\in X_U:=B^\times\bs\wh{B}^\times/U$, let $L_h^o\subset L_h:=h\wh{L}h^{-1}\cap V$ be the subset
\begin{equation}
L_h^o=\left\{\ell\in L_h\ \Big|\ \ell\in h\cdot L_p^o\text{ for all $p|2M$}\right\}.\label{ort}\tag{ot}
\end{equation}
 We refer to $L_h^o\subset L_h$ as the oriented solutions of the underlying ternary quadratic form.
\end{defn}

In light of  \eqref{FP}, \eqref{FT}, it now follows: For each $D\in\Sigma_C $, normalized Fourier coefficients of $\theta_C$ satisfy
\[a_{|D|}(\theta_C)=\frac{P_{K_D}^0(f_D)}{[\CO_D^\times:\BZ^\times]}\] and
\[a_{|D|}(\theta_C)=\sum_{[h]\in X_U}\frac{f(h)}{w_h}\cdot\#  (L_h^o\cap V_{D}(\BQ)),\]
where $w_h:=\#(hUh^{-1}\cap H(\BQ))$ is a finite group.

\subsubsection{Main result}

Let $\pi\subset  \CA_{0,2}(H)$ and $M$ the conductor of the Jacquet--Langlands transfer. Let $C$ be an  $M$-equivalence class with $\varepsilon(C)=1$ such that $(\pi, C)$ satisfies Tunnell--Saito condition \eqref{TS}.

 \begin{thm}\label{thmm}Assume $B$ is definite and $L(\frac{1}{2},\pi)\neq 0$.  Let $(f,\phi)$ be a distinguished pair as in \S \ref{distp}.
 
 Then there exists an explicit constant $k_C\in \BQ^\times$ such that for each $D\in \Sigma_C$,
\[\frac{ \sqrt{D}\cdot L\left(\frac{1}{2},\pi\otimes\chi_{D}\right)}{ \Omega_\sigma^-}=k_C\cdot  \left(\sum_{[h]\in X_U}\frac{f(h)}{w_h}\cdot\#  (L_h^o\cap V_{D}(\BQ))\right)^2.\] Here the notation is as in \S \ref{distp}. (cf. \cite{HTX})
\end{thm}
\begin{remark}\label{rem1}\
\begin{itemize}
\item[(i)]
 A more general set-up appears in \cite{HTX}: weight $k$ and indefinite quaternion algebras\footnote{Accordingly, the infinite component $\pi_\infty$ may be non-trivial and in general, one needs to account for infinitely many integral solutions with fixed
reduced norm in an indefinite lattice.}.

\item[(ii)]

 If for each $v|M$, $\varepsilon(\pi_v)=\epsilon(B_v)$, then the local obstruction of Atkin--Lehner operator disappears. The oriented solutions in Theorem \ref{thmm} may thus be replaced by the whole lattice solutions - switching the constant $k_C$ with a $2$-power multiple.

 \item[(iii)]
     If $\CE$ is CM or there exists $A\in\CE$ with non-square conductor, then there exists $n_1>0$, $n_2<0$: both $L(s,E^{(n_1)})$ and $L(s,E^{(n_2)})$ have sign $+1$ and also do not vanish at the center. In this case, there exists a partition of $\CE$ such that each $\CE_C$ corresponds to a definite quaternion algebra.
\end{itemize}
\end{remark}

\subsection{Congruent Number $L$-values, revisited }As an application of Theorem \ref{thmm}, we recover  Tunnell's theorem and exhibit a new interpretation of the central $L$-values of the congruent number elliptic curves
in terms of ternary quadratic forms.

\subsubsection{Tunnell's theorem, again}\label{ss:Tuag} Let $\sigma$ be the irreducible cuspdial automorphhic representation associated to $E:y^2=x^3-x$, its conductor $M=32$.
We first use Theorem \ref{walds} to recover Tunnell's result.
  Let $C\subset\BQ_{<0}$ be an $M$-equivalence class consisting of negative rationals. Then $\varepsilon(C)=1$ if and only if $C=[-1], [-2], [-10],  [-3]$, where $[n]$ denotes the $M$-equivalence class of $n\in \BQ^\times$.
Now pick $C$ with $\varepsilon(C)=1$ to be one of $[-1], [-2], [-3]$.

The quaternion algebra determined by $C$ such that $(\sigma,C)$ satisfies Tunnell--Saito condition \eqref{TS} is the quaternion algebra $B$ over $\BQ$ ramified exactly at $2$ and infinity, i.e.
\[B=\BQ\oplus \BQ i\oplus  \BQ j\oplus \BQ k\] with $i^2=j^2=-1$, $ij=k=-ji$. Let $H=PB^\times$ and $\pi\subset \CA_{0,2}(H)$ be irreducible: the Jacquet--Langlands correspondence of $\pi$ is $\sigma$. For the existence of admissible pairs, we begin with the computation of  the space $V(\pi,  K)$ (cf.~\cite{TYZ}).

Let $K=\BQ(x)\subseteq B$ with
$$x=i,\quad  i+j,\quad  i+j+k, \qquad\text{for}\ C=[-1],\  [-2],\   [-3], \quad \text{respectively}.$$
Let $\CO_B$ be the maximal order of $B$ generated over $\BZ$ by $ i, j, \displaystyle{\frac{1+i+j+k}{2}}$. Then
$$R=\CO_K+4\CO_B$$
is an admissible order for $(\pi, {\bf 1}_K)$.  The (finite) Shimura set $X_{\wh{R}^\times}$ of level $\wh{R}^\times$ has representatives given by elements in $H(\BQ_2):$
\[X_{\wh{R}^\times}=\begin{cases}
\{1_2, (1+2j)_2, (1+1+i+j+k)_2, (1+1+i-j+k)_2\},& \text{if $C=[-1]$,}\\
\{1_2, (1+2j)_2\},& \text{if $C=[-2]$,}\\
\{1_2, (1+2j)_2\},& \text{if $C=[-3]$.}
 \end{cases}\]
 
 A basis $f: X_{\wh{R}^\times} \ra \BZ$ of the one dimensional space $V(\pi,K)$
 with respect to the above representatives of $X_{\wh{R}^\times}$ (cf.~\cite{TYZ}):
\[f=\begin{cases}
(1,-1,0,0),& \text{if $C=[-1]$,}\\
(1,-1),& \text{if $C=[-2]$,}\\
(3,-1),& \text{if $C=[-3]$.}
 \end{cases}
 \]
Let $L\subset V$ be the lattice with level $\begin{cases}
  128,&\quad \text{if $C=[-1], [-3]$}\\
  256,&\quad \text{if $C=[-2]$}\\
\end{cases}$ and character $\mathbf{1}$ given by\[L=\begin{cases}
 \BZ i\oplus \BZ(j+k)\oplus \BZ4(j-k),&\quad \text{if $C=[-1], [-3]$}\\
 \BZ (i+j)\oplus \BZ 2(i-j)\oplus \BZ8k,&\quad \text{if $C=[-2]$.}
\end{cases}\] Put $\phi=\textbf{1}_{\wh{L}}\otimes e^{2\pi |c| q(\cdot)}\in \CS(V(\BA))$. Then $(f,\phi)$ is an admissible pair.  

In view of Theorem \ref{walds} and \eqref{FT}, finally Tunnell's result:

 For any positive square-free integer $n$,
      \[\frac{L(1,E^{(n)})}{\Omega/\sqrt{n}}=\frac{1}{16}\cdot\begin{cases}
      (\# L_1\cap V_{-n}(\BQ)-\#L_2\cap V_{-n}(\BQ))^2,&\quad \text{if $n\equiv 1,3\pmod{8}$}\\
      2(\# L_1'\cap V_{-n}(\BQ)-\#L_2'\cap V_{-n}(\BQ))^2,&\quad \text{if $n\equiv 2\pmod{8}$.}
      \end{cases}
      \] Here  $\Omega=\int_1^\infty\frac{dx}{\sqrt{x^3-x}}$,
$\{L_1,L_2\}$, $\{L_1',L_2'\}$ the genus class of lattices in $V$ given by
\[\{L_1=\BZ i\oplus \BZ(j+k)\oplus \BZ4(j-k), L_2=\BZ(i+j)\oplus\BZ 2k\oplus\BZ (2i-2j+k)\}\] and
\[\{L_1'=\BZ (i+j)\oplus \BZ 2(i-j)\oplus \BZ8k,L_2'=\BZ2(i-j)\oplus \BZ 2(i+j)\oplus\BZ(i+j-4k)\}.\]
Note that the ternary quadratic forms corresponding to $L_1,L_2,L_1',L_2'$ are given by
\[x^2+2y^2+32z^2,\quad 2x^2+4y^2+9z^2+4yz,\]
\[2x^2+8y^2+64z^2,\quad 8x^2+8y^2+18z^2+8yz,\]respectively.

\subsubsection{Congruent number L-values, anew} In light of  Theorem \ref{thmm}, one may obtain central $L$-value formulae for congruent number elliptic curves in term of different ternary quadratic forms (or equivalently different lattices).
\begin{itemize}
\item[(I)]\label{A}
 Let $x=2i\in B$ and $R=\BZ[i]+4\CO_B$. Let $f$ be the test vector in the last subsection for $C=[-1]$. Consider  the genus class of lattices in $V$ with level $128$ and character $\mathbf{1}$ consisting of
 $$\begin{aligned}
&L_1=\BZ i\oplus \BZ4(j-k)\oplus\BZ4(j+k),\qquad \\
&L_2=\BZ2i\oplus\BZ(i+4k)\oplus\BZ(4j-i), \\
&L_3=\BZ2i\oplus\BZ(i+2j+2k)\oplus \BZ4(j-k).\end{aligned}$$
  Then the lattice $2L_1$ is $x$-distinguished and $\displaystyle(f,\textbf{1}_{2\wh{L}_1}\otimes e^{2\pi |c| q(\cdot)})$ a distinguished pair.
\begin{prop}\label{propA}
For any positive square free integer $n\equiv 1\pmod{8}$,
      \[\frac{L(1,E^{(n)})}{\Omega/\sqrt{n}}=\frac{1}{16}\big(\# L_1\cap V_{-n}(\BQ)-\#L_2\cap V_{-n}(\BQ)\big)^2.\]
\end{prop}

\item[(II)]\label{B}We now construct a new distinguished pair whose Fourier coefficients interpret the 
$L$-values $L(1,E^{(n)})$ for positive square-free integers $n\equiv 1\pmod{8}$.

Let $\sigma'=\sigma\otimes\chi_2$ be with conductor $M'=64$. Consider the $M'$-equivalence class $C=[-2]$ with $\varepsilon(C)=1$. The pair $(\sigma',C)$ gives the quaternion algebra $B$
as in \S\ref{ss:Tuag}.
Let $\pi':=\pi\otimes\chi_2\subset \CA_{0,2}(H)$ and $K:=\BQ(x)\subset B$ with $x:=2(i+j)$.

Then
$$R:=\CO_K+4(i+j)\CO_B$$
is an admissible order for $(\pi', {\bf 1}_K)$. The Shimura set $X_{\wh{R}^\times}$ has representatives given by elements in $H(\BQ_2):$
\[X_{\wh{R}^\times}=\left\{1_2,(1+2j)_2,(1+2(1+i+j+k))_2,(1+3(1+i+j+k))_2\right\}.\]
A $\BZ$-primitive test vector in $V(\pi',K)$ is given by
\[f'=(1,1,-1,-1).\]
Consider  the genus class of lattices in $V$ with level $512$ and character $\mathbf{1}$ consisting of
 $$\begin{aligned}
&L_1'=\BZ (i+j)\oplus \BZ8(i-j)\oplus\BZ8k,\qquad \\
&L_2'=\BZ2(i+j)\oplus\BZ8k\oplus\BZ(-3i+5j-4k), \\
&L_3'=\BZ8k\oplus\BZ16i\oplus \BZ(9i+j),\\
&L_4'=\BZ2(i+j)\oplus\BZ8(i-j)\oplus \BZ(i+j+4k).
\end{aligned}$$
Then the lattice $2L_1'$ is $x$-distinguished and $\displaystyle(f',\textbf{1}_{2\wh{L_1'}}\otimes e^{2\pi |c| q(\cdot)})$ a distinguished pair.
\begin{prop}\label{propB}
   For any positive square-free integer $n\equiv 1\pmod{8}$,
      \[\frac{L(1,E^{(n)})}{\Omega/\sqrt{n}}=\frac{1}{16}\big(\# L_1'\cap V_{-2n}(\BQ)+\# L_2'\cap V_{-2n}(\BQ)-\# L_3'\cap V_{-2n}(\BQ)-\#L_4'\cap V_{-2n}(\BQ)\big)^2.\]
\end{prop}
\end{itemize}
\begin{remark}\label{nw} Let $\theta$, $\theta'$ be the theta lifting of the distinguished pairs as in (I), (II) respectively.
For positive square-free $n\equiv 1\pmod{8}$, one has 
\[a_{4n}^2(\theta)=a_{8n}^2(\theta')\] 
(cf. Proposition \ref{propA} and Proposition \ref{propB}).
However, as $n$ varies in positive square-free integers with $n\equiv 1\pmod{8}$, note that the sequence 
$\{a_{4n}(\theta)\}$ is not proportional\footnote{For example, $a_{4}(\theta)=a_{8}(\theta')$ but $a_{4\cdot57}(\theta)=-a_{8\cdot57}(\theta')$.} to $\{a_{8n}(\theta')\}$. 
\end{remark}
\section{$p$-converse}\label{s:pcv}\subsubsection{The Birch and Swinnerton-Dyer conjecture, bis}
Let $A$ be an elliptic curve over $\BQ$ and $p$ a prime.

The $p^{\infty}$-Selmer group $\Sel_{p^{\infty}}(A/\BQ)$
appears as the middle term of the short exact sequence
\begin{equation}\label{ses}
0 \ra A(\BQ) \otimes_{\BZ} \BQ_{p}/\BZ_{p} \ra \Sel_{p^{\infty}}(A/\BQ) \ra \Sha(A/\BQ)[p^{\infty}] \ra 0,
\end{equation}
for $\Sha(A/\BQ)[p^{\infty}]$ the $p$-primary part of $\Sha(A/\BQ)$.

In view of the exact sequence \eqref{ses}, Conjecture \ref{cBSD} suggests:

\begin{conj}\label{sBSD}
Let $A$ be an elliptic curve over $\BQ$. The following are equivalent. 
\begin{itemize}
\item[(a)] $\rank_{\BZ}A(\BQ)=r$ and $\Sha(A/\BQ)$ is finite.
\item[(b)] $\corank_{\BZ_{p}} \Sel_{p^{\infty}}(A/\BQ)=r$ for $p$ a prime.
\item[(c)] $\ord_{s=1}L(s,A)=r$.
\end{itemize}
\end{conj}

Part (b) follows from part (a) just by \eqref{ses}. That $(c) \implies (a)$ is a spectacular result\footnote{For a brief introduction, one may see \cite[\S3.1--3.2]{BSTsv}.}
 towards the Birch and Swinnerton-Dyer conjecture due to Coates--Wiles \cite{CoWi} and Rubin \cite{Ru0} (the CM case), and
Gross--Zagier \cite{GZ} and Kolyvagin \cite{Ko} (the general case).

After Skinner, nowadays, `$(b) \implies (c)$' is referred to as a $p$-converse: a $p$-adic criterion to have analytic rank $r$. For $r=0$, an important progress towards the $p$-converse - Rubin \cite{Ru} (the CM case) and Skinner--Urban \cite{SU} (the ordinary non-CM case). The $r=1$ case remained widely open until the breakthrough due to  Zhang \cite{Zh} and Skinner \cite{Sk} (the ordinary non-CM case) a few years back. Since then, the $p$-converse is undergoing a revival, in light of which one may hope a complete resolution of:
\begin{conj}[$p$-converse]\label{pcv} Let $A$ be an elliptic curve over $\BQ$, $p$ a prime and $r=0, 1$. Then,
$$\corank_{\BZ_p} \Sel_{p^\infty} (A/\BQ)=r \quad \Longrightarrow\quad \ord_{s=1} L(s,A)=r.$$
	\end{conj}
Our study concerns some of the missing cases, notably the CM curves (cf. \cite{TA1}, \cite{BuTi2}, \cite{BST1}, \cite{BST2}, \cite{BST4}). For a brief overview of the current progress towards Conjecture \ref{pcv}, one may refer to \cite[\S3.1--3.2]{BSTsv}.
\subsubsection{The Goldfeld conjecture, bis}
In view of Conjecture \ref{Gd} and Conjecture \ref{pcv}:
\begin{conj}\label{GdSel}Let $A$ be an elliptic curve over $\BQ$ and $p$ a prime.

Then, for a density one subset of square-free integers $d $ with $\varepsilon(A^{(d)})=+1$ (resp.~$\varepsilon(A^{(d)})=-1):$
\[
\corank_{\BZ_{p}}\Sel_{p^{\infty}}(A^{(d)}/\BQ)=0
,\quad \text{(resp.~$\corank_{\BZ_{p}}\Sel_{p^{\infty}}(A^{(d)}/\BQ)=1
$)}.\]

\end{conj}
\subsubsection{$p$-converse to a theorem of Coates--Wiles and Rubin}

\begin{thm}\label{pcCWR}
Let $A$ be a CM elliptic curve over $\BQ$ and $p$ a prime.
Then
$$
\mathrm{corank}_{\BZ_p}\Sel_{p^{\infty}}(A/\BQ) =0  \implies L(1,A) \neq 0.
$$
\end{thm}

In the early 1990's Rubin \cite{Ru} proved this $p$-converse when $p \nmid \# \CO_{K}^{\times}$ for $K$ the CM field -
the hypothesis being essential to utilise the Euler system of elliptic units. The case $p=2$  remained open since then. The unconditional $p$-converse is quite recent \cite{BuTi2}.

Now, in view of Theorem \ref{Tun}:
\begin{cor}\label{Qeq}
For a positive square-free integer $n\equiv 1,2,3\mod{8}$, let $E^{(n)}$ be the congruent elliptic curve $ny^2=x^3-x$. Let $p$ be a prime. Then, 
$$
\CL(n)\neq 0 \iff \mathrm{corank}_{\BZ_p}\Sel_{p^{\infty}}(E^{(n)}/\BQ) =0.
$$
\end{cor}
\vskip2mm
{\it{Approach}.} Unconventionally - for the CM case - the approach employs the Beilinson--Kato elements \cite{K}.
It is inherently Iwasawa theoretic, principle:
$$
\text{Galois actions on arithmetic objects $\longleftrightarrow$ zeta values.}
$$
(cf. \cite{K1}).
Indeed, the decisive shift is to consider Kato's main conjecture for $A$ (cf.~\cite[Conj.~12.10]{K}) instead of the habitual elliptic units main conjecture for $K$ (cf. \cite{Ru}). The heart of \cite{BuTi2} is the proof of Kato's main conjecture, which via an Iwasawa descent implies Theorem \ref{pcCWR}.

\begin{remark}
Unlike \cite{Ru}, our approach uniformly treats the ordinary and non-ordinary primes.
\end{remark}
\subsubsection{$p$-converse to a theorem of Gross--Zagier, Kolyvagin and Rubin}
\begin{thm}\label{pcGZKIII}
Let $A$ be a CM elliptic curve over the rationals with conductor $N$ and $p\nmid 6N$ a prime. Then,
$$
\corank_{\BZ_{p}}\Sel_{p^{\infty}}(A/\BQ)=1 \implies \ord_{s=1}L(s, A)=1.
$$
\end{thm}

For $p$ also a prime of ordinary reduction, this $p$-converse was proved
 in \cite{BuTi2}. It is a rare instance where the non-CM case \cite{Sk}, \cite{Zh} preceded the CM case.
 The supersingular case will appear in \cite{BST4} (cf.~\cite{BKO1}, \cite{BKO2}).
 Another approach, which generalizes to CM curves over totally real fields, is given in \cite{BST2}.

 A salient feature of \cite{BuTi2}, \cite{BST2}: the arithmetic of auxiliary Heegner points over the CM field. For CM curves, the CM field does not satisfy the Heegner hypothesis - while - it is natural to seek to utilise the arithmetic of the CM field.  The generality of the Gross--Zagier formula \cite{YZZ} allows us to still introduce auxiliary Heegner points over the CM field.
 The core of the approach: an interplay between Iwasawa theory of the auxiliary Heegner points \cite{PR2} and CM Iwasawa theory \cite{Ru}. For a more detailed account of the strategy, one may refer to \cite[\S4]{BSTsv}.

\section{Distribution of Selmer groups}\label{s:dsel}

\subsection{Conjectures}
The subsection presents conjectures on the distribution of Selmer groups of elliptic curves over a fixed number field following \cite{PR}, \cite{BKLPR}. As an instructive introduction, one may refer to \cite{P}.

Let $A$ be an elliptic curve over a number field $F$ and $p$ a prime. For the $p^{\infty}$-Selmer group $\Sel_{p^\infty}(A/F)$, the exact sequence
\begin{equation}\label{$Seq_{A}$}\tag{$\Seq_{A}$}
 0\ra A(F)\otimes \BQ_p/\BZ_p\xrightarrow\kappa \Sel_{p^\infty}(A/F)\ra \Sha(A/F)[p^\infty]\ra 0
\end{equation}
 of cofinitely generated $\BZ_p$-modules is split.
 \subsubsection{The mod $p$ Selmer groups}\label{ss:pSel}
We recall the distribution model for the mod $p$ Selmer groups $\Sel_{p}(A/F)$ due to Poonen and Rains \cite{PR}.

Consider the infinite dimensional locally compact hyperbolic quadratic space
  $$V:={\prod_v}' H^1(F_v, A[p]).$$ Here the restricted product arises from $\{A(F_v)/pA(F_v)\}_v$ and the quadratic form $Q:$ the sum of local quadratic forms $Q_v$ corresponding to the map
  \[H^1(F_v,A[p])\ra H^2(F_v,\BG_m)\hookrightarrow \BQ/\BZ,\]
  which
  arises from the short exact sequence
  \[0\ra \BG_m\ra \CH\ra A[p]\ra 0\]
 of $G_{F_v}$-modules for $\CH$ the Heisenberg group scheme as in \cite{PR}.

By definition, $\Sel_p(A/F)$ is the intersection of two maximal isotropic subspaces:
 the images of $H^1(F, A[p])$ under the restriction map and the local Kummer maps:
  \[\xymatrix{
  \prod_{v}A(F_v)/pA(F_v)\ar[dr]^{\kappa_{v}}&&\\
  &\prod_v' H^1(F_v, A[p])&\\
  H^1(F,A[p])\ar[ur]^{\Res}&&
  }\]This perhaps suggests the following model.

  Let
  $(V=W\oplus W^\vee,Q)$ be a hyperbolic quadratic space over $\BF_p$ of dimension $2n$ with the natural quadratic form and $I_V$ the set of maximal isotropic subspaces.

  A key proposal \cite{PR}:

  \begin{conj}\label{PR}
  \[\Prob\left(\dim \Sel_p(A/F)=d\right)=\lim_{\dim V\ra \infty}\Prob\left(\dim (Z_1\cap Z_2)=d\ \Big|\ Z_1, Z_2\in I_V\right)\]

  \end{conj}It may be seen that
    \[C_{p,d}:=\lim_{\dim V\to\infty}\Prob\left(\dim (Z_1\cap Z_2)=d\ \Big|\ Z_1, Z_2\in I_V\right)=\prod_{j=0}^\infty(1+p^{-j})^{-1}\prod_{i=1}^d \frac{p}{p^i-1}.\]
In particular, when ordered by height, the elliptic curves with $\rank_{\BZ}A(F)\geq 2$ have density $\leq \frac{p+1}{p^2}$.
In light of the BSD conjecture, the Poonen--Rains conjecture thus implies the following.
  \begin{conj}[Rank conjecture]\label{rank}Let $r\in\{0,1\}$. When ordered by height, $50\%$ of the elliptic curves over $F$ have Mordell--Weil rank $r$.
  \end{conj}

\subsubsection{The $p^{\infty}$-Selmer groups}
 Following Bhargava, Kane, Lenstra, Poonen and Rains (cf. \cite{BKLPR}), we now switch to the $p^{\infty}$-Selmer group.

Equip $\BZ_p^{2n}$ with the hyperbolic quadratic form $Q:\BZ_p^{2n}\ra \BZ_p$ given by
\[Q(x_1,\cdots,x_n,y_1,\cdots,y_n)=\sum_{i=1}^nx_i y_i.\] A direct summand $Z$ of $\BZ_p^{2n}$ is called maximal isotropic if $Q|_Z=0$ and the rank of $Z$ is $n$.
Let $\OGr_n(\BZ_p)$ be the set of the maximal isotropic $\BZ_p$-submodules of $\BZ_p^{2n}$.
Let $\OGr_n$ be the underlying smooth projective scheme over $\BZ$. Consider the natural probability measure $\nu_n$ on $\OGr_{n}(\BZ_p):$
\[\nu_n(S):=\lim_{e\ra\infty}\frac{\#\Im\left(S\ra \OGr_{n}(\BZ/p^e\BZ)\right)}{\#\OGr_{n}(\BZ/p^e\BZ)}\]
for $S\subset \OGr_n(\BZ_p)$ an open and closed subset. Fix $W:=\BZ_p^n\times 0$ and $Z$ at random\footnote{In fact, equivalent to choose both $Z$ and $W$ at random, since $O(V, Q)$ acts transitively on $\OGr_n(\BZ_p)$.}.

  In the spirit of \S\ref{ss:pSel} define
\[\begin{aligned}
&R:=(Z\cap W)\otimes_{\BZ_p}\BQ_p/\BZ_p,\\
&S:=(Z\otimes_{\BZ_p}\BQ_p/\BZ_p)\cap (W\otimes_{\BZ_p}\BQ_p/\BZ_p), \\
&T:=S/R.\end{aligned}\]
  Note that $T$ is finite, endowed with a canonical non-degenerate alternate pairing. Let $\CQ_{2n}$ be the distribution of the isomorphism class of the short exact sequences
  \begin{equation}\label{rses}
  0\ra R \ra S\ra T\ra 0,
  \end{equation}
induced from $\nu_n$.
It may be seen that the limit of $\CQ_{2n}$ exists, say $\CQ$.

A fundamental proposal \cite{BKLPR}:

\begin{conj}\label{cBKLPR}
Let $F$ be a global field and $\CS$ a short exact sequence of $\BZ_p$-modules as in \eqref{rses}. When  ordered by height,
the density of \[\big{\{}A:\Seq_{A}\simeq \CS\big{\}}\] equals the $\CQ$-probability of $\CS$.

More precisely, let $G$ be a finite symplectic $p$-group.
Then
\[
\operatorname{Prob}\Big( \Sel_{p^\infty}(A/F)\simeq (\BQ_p/\BZ_p)^r\oplus G\Big)\\
=\frac{(\#G)^{1-r}}{\#\Sp(G)}\times\begin{cases}
\displaystyle{\frac12\prod_{i=r+1}^\infty(1-p^{1-2i})},
&\text{if }r=0,1, \\
0,&\text{if }r\geq 2.
\end{cases}\]
\end{conj}

The scheme $\OGr_{n}$ is a disjoint union of two isomorphic subschemes $\OGr_n^\pm$ such that - for any  field $k$ - $\OGr_n^+(k)$ (resp.~$\OGr_n^-(k)$) parameterize $Z\in \OGr_n(k)$ with $\dim (Z\cap W_k)$ even (resp.~odd). Further, the locus with $\dim (Z\cap W_k)\geq 2$ has lower dimension (cf. \cite{BKLPR}). Thus Conjecture \ref{cBKLPR} implies Conjecture \ref{rank}. Assuming independence in $p$, it  also implies:

\begin{conj}\label{conj1.3}

Let $n\in \BZ_{\geq 1}$. When ordered by height, the average of $\#\Sel_n(A/F)$ is $\sigma(n):=\sum_{d|n} d$.
\end{conj}
In particular, the average of $\#\Sel_{2^n}(A/F)$ is conjecturally $2^{n+1}-1$, which resonates through the following subsection.

\subsection{Smith's work}
\label{s:smith}
In his remarkable thesis, Smith \cite{Sm1} made the following breakthrough towards the distribution of the
$2^{\infty}$-Selmer groups in the quadratic twist family of elliptic curves over $\BQ$ (cf. \cite{GM}).
\begin{thm}\label{Sm}
Let $A/\BQ$ be an elliptic curve such that
\begin{equation}\label{ord}\tag{ord}
\text{$A[2]\subset A(\BQ)$ and $A$ has no cyclic subgroup of order $4$ defined over $\BQ$.}
\end{equation}
Let $r\in \BZ_{\geq 0}$ and $G$ be a symplectic $2$-group.

Then, as $d$ varies over square-free integers:
\[
\operatorname{Prob}\left(
\Sel_{2^\infty}(A^{(d)}/\BQ)
\cong(\BQ_2/\BZ_2)^r\oplus G\right)
=\frac{(\#G)^{1-r}}{\#\Sp(G)}\times\begin{cases}
\frac12\prod_{i=r+1}^\infty(1-2^{1-2i}),
&\text{if }r=0,1, \\
0,&\text{if }r\geq 2.
\end{cases}
\]

\end{thm}

A couple of striking consequences:

\begin{cor}\label{SmGo}Let $A$ be an elliptic curve over $\BQ$ as in Theorem \ref{Sm}.

Then, for a density one subset of square-free integers $d$ with $\varepsilon(A^{(d)})=+1$ (resp.~$\varepsilon(A^{(d)})=-1):$
\[
\corank_{\BZ_{2}}\Sel_{2^{\infty}}(A^{(d)}/\BQ)=0
,\quad \text{(resp.~$\corank_{\BZ_{2}}\Sel_{2^{\infty}}(A^{(d)}/\BQ)=1
$)}.\]

\end{cor}

\begin{cor}\label{SmCo}
The density of non-congruent numbers in all positive square-free integers $n\equiv 1,2,3 \mod{8}$ is one.
\end{cor}
\begin{remark}
The quadratic twist family of elliptic curves is complementary to Conjecture \ref{cBKLPR}, yet intriguingly Theorem \ref{Sm} echoes the same principal.
\end{remark}
\begin{remark}\label{SmFo}\noindent
\begin{itemize}
\item[$\circ$] Very recently, Smith has announced: Theorem \ref{Sm} also holds for elliptic curves $A/\BQ$ with $A(\BQ)[2]=0$. When ordered by height, a density one subset of elliptic curves over $\BQ$ satisfies the hypothesis.
\item[$\circ$] Suppose $A(\BQ)[2]\simeq \BZ/2\BZ$ and let $A'$ be the isogenous curve arising from the $2$-torsion. If $\BQ(A[2])\neq \BQ(A'[2])$, Smith has announced: Corollary \ref{SmGo} still holds.
\end{itemize}
\end{remark}

\subsection{Goldfeld's conjecture: an instance}
By Corollary \ref{SmGo} and Theorem \ref{pcCWR},

\begin{cor}\label{eGd}
Let $E^{(n)}:y^{2}=x^{3}-n^{2}x, n\in \BN,$ be  congruent number elliptic curves. Then,
$$
\text{$L(1,E^{(n)})\neq 0$ for a density one subset of positive square-free integers $n \equiv 1,2,3 \mod{8}$.}
$$
\end{cor}
\begin{remark}
In light of Remark \ref{SmFo} and Theorem \ref{pcCWR}: The even parity Goldfeld conjecture holds for quadratic twist family of CM elliptic curves over $\BQ$, once the CM field differs from $\BQ(\sqrt{-2})$.
\end{remark}

\section{Distribution of $2$-Selmer groups: exceptional case}\label{s:dex}

Let $E$ be a fixed elliptic curve over $\BQ$. Recall that the $2$-Selmer group
$\Sel_2(E/\BQ)$
appears in the fundamental short exact sequence
$$
0\to E(\BQ)/2E(\BQ)\xrightarrow\kappa\Sel_2(E/\BQ)\to\Sha(E/\BQ)[2]\to 0.
$$
In particular, if $\Sel_2(E/\BQ)=\kappa(E(\BQ)_\tor)$,
then $E$ has Mordell--Weil rank $0$ and $\Sha(E/\BQ)[2]=0$.
On the other hand,
if $\Sel_2(E/\BQ)/\kappa(E(\BQ)_\tor)$ is non-trivial and $E$ has Mordell--Weil rank $0:$
$\Sha(E/\BQ)[2]$ is non-trivial.

The $2$-parity conjecture, proved  by Monsky \cite{Mon96} (see also \cite{DD1}), asserts that $$\dim_{\BF_2}\big(\Sel_2(E/\BQ)/\kappa(E(\BQ)_\tor)\big)\equiv
\displaystyle\operatorname*{ord}_{s=1}L(s,E)\mod 2.$$

 For simplicity assume that
$E[2]\subset E(\BQ)$. For a square-free integer $n$, let $E^{(n)}$ be the quadratic twist of $E$ by $\BQ(\sqrt{n})$ and define the rational number \[\displaystyle \mathscr{L}(n):=\frac{L(1,E^{(n)})}{\Omega_{E^{(n)}}}\cdot \left(\frac{\prod_{\ell}c_\ell(E^{(n)})}{\# E^{(n)}(\BQ)_{\tor}^2}\right)^{-1},\] 
the analytic Sha of $E^{(n)}$ if $L(1,E^{(n)})\neq 0$.

Given a residue class $a\mod M$, one may seek to study: among
all positive (or negative) square-free integers $n$ in this residue class,
the minimal value of $\dim_{\BF_2}\Sel_2(E^{(n)}/\BQ)/\kappa(E^{(n)}(\BQ)_\tor)$
and the $2$-adic valuation of $\mathscr{L}(n)$ of $E^{(n)}$,
and the distribution of $n$ which $\ord_2 \mathscr{L}(n)$
attains the minimum.
\subsubsection{A precursor}
If $E$ satisfies \eqref{ord},
then the results of
Heath-Brown \cite{HB2}, Swinnerton-Dyer \cite{SD} and Kane \cite{Kane}
exhibit the minimal value of $2$-Selmer ranks
and its distribution among the residue classes.
The results suggest that
the distribution of $2$-Selmer groups in a quadratic twist family
(modulo the contribution of torsion points) still mirrors the Poonen--Rains principle\footnote{though their framework excludes the quadratic twist family} (cf. Conjecture \ref{PR}).
The phenomenon inspired and resonates through
Smith's work on $2^\infty$-Selmer groups \cite{Sm2}, \cite{Sm1}.

The main results of \cite{HB2}, \cite{SD}, \cite{Kane}:

\begin{thm}
\label{2-Sel-ordinary}
Let $E/\BQ$ be an elliptic curve
satisfying \eqref{ord}.
Let $n_0$ be a square-free integer and
$d_0=\dim_{\BF_2}
\Sel_2(E^{(n_0)}/\BQ)/\kappa(E^{(n_0)}(\BQ)_\tor)$.
Let $N$ be the conductor of $E$,
and $[n_0]\subset\BQ^\times$ the $N$-equivalence class\footnote{See Definition \ref{M-equivalent}.}
which contains $n_0$.

{\rm(i)}
For $d,k \in \BZ_{\geq 0}$,
 let $\pi_{d,k}$ denote
the following probability
among square-free integers $n$
with $k$ distinct prime factors
such that
$n\in[n_0]:$
$$
\pi_{d,k}:=\operatorname{Prob}\left(\dim_{\BF_2}
\Sel_2(E^{(n)}/\BQ)/\kappa(E^{(n)}(\BQ)_\tor)=d\right).
$$
Then
$$
\lim_{k\to\infty}\pi_{d,k}
=\begin{cases}
2\prod_{j=0}^\infty(1+2^{-j})^{-1}\prod_{i=1}^d \frac{2}{2^i-1},&\text{if }
d\equiv d_0\pmod 2, \\
0,&\text{if }d\not\equiv d_0\pmod 2.
\end{cases}
$$

{\rm(ii)}
For $d\in \BZ_{\geq 0}$,
among square-free integers $n\in[n_0]:$
$$
\operatorname{Prob}\left(\dim_{\BF_2}
\Sel_2(E^{(n)}/\BQ)/\kappa(E^{(n)}(\BQ)_\tor)=d\right)
=\begin{cases}
2\prod_{j=0}^\infty(1+2^{-j})^{-1}\prod_{i=1}^d \frac{2}{2^i-1},&\text{if }
d\equiv d_0\pmod 2, \\
0,&\text{if }d\not\equiv d_0\pmod 2.
\end{cases}
$$
In particular,
for $d\in \BZ_{\geq 0}$,
among square-free integers $n:$
$$
\operatorname{Prob}\left(\dim_{\BF_2}
\Sel_2(E^{(n)}/\BQ)/\kappa(E^{(n)}(\BQ)_\tor)=d\right)
=\prod_{j=0}^\infty(1+2^{-j})^{-1}\prod_{i=1}^d \frac{2}{2^i-1}.
$$
\end{thm}

Part (i) is due to Swinnerton-Dyer \cite{SD},
a key: $(\pi_{d,k})_{d=0}^\infty$ is connected by a Markov chain
as $k$ varies.
Via analytic tools, Kane \cite{Kane} proved that part (i) implies part (ii)
- a transition from the density for integers
with restricted prime factors to natural density.
The above theorem for the congruent number elliptic curve $y^2=x^3-x$
is also due to Heath-Brown \cite{HB2}, an independent approach.
\subsubsection{An exceptional case}
Some quadratic twist families of elliptic curves over $\BQ$ satisfy neither \eqref{ord}
nor the hypotheses in Remark \ref{SmFo}.

A key missing case:
\begin{equation}\label{exc}\tag{exc}
\text{$E[2]\subset E(\BQ)$ and $E$ has rational $4$-torsion points.}
\end{equation}
For example, the quadratic twist family of tiling number elliptic curves 
\[E^{(n)}:y^2=x(x-n)(x+3n).\]
Notice $E(\BQ)_\tor\cong\BZ/2\BZ\times\BZ/4\BZ$.

Perhaps surprisingly, in light of the presence of such rational $4$-torsion, the distribution
of $2$-Selmer groups no longer seems to be as in Theorem \ref{2-Sel-ordinary}.
A suggestive example \cite{FLPT}:
\begin{prop}
Let $E$ be the elliptic curve $y^2=x(x-1)(x+3)$.
If $n\neq 1$, $n\equiv 1\pmod{12}$ is positive square-free,
then
$$
\dim_{\BF_2}
\Sel_2(E^{(-n)}/\BQ)/\kappa(E^{(-n)}(\BQ)_\tor)\geq 2,
$$
and the corresponding $\mathscr{L}(-n)$ is also even.
\end{prop}
\subsection{Main results}
Our preliminary study suggests that for elliptic curves satisfying \eqref{exc},
the distribution of $2$-Selmer groups may resemble that of the
$4$-ranks of ideal class groups of the underlying imaginary quadratic fields.

Let $g(n):=\#2\Cl(\BQ(\sqrt{-n}))$,
in particular, $g(n)$ is odd if and only if $\BQ(\sqrt{-n})$
has no ideal class of exact order $4$.

A main result of \cite{FLPT}:

\begin{thm}
\label{FLPT-TFAE}
Let $E$ be the elliptic curve $y^2=x(x-1)(x+3)$.
Let $n\equiv 3,7\pmod{24}$ be a positive square-free integer.
Let $\epsilon\in\{\pm 1\}$.
Then the followings are equivalent. 
\begin{enumerate}
\item[(a)]
The genus invariant
$$
\begin{cases}\displaystyle
g(n)+\sum_{\substack{d\mid n\\ d\equiv 11\mod 24}}g(n/d)g(d),&\text{if
$n\equiv 7\mod 24$ and $\epsilon=-1$}, \\
g(n),&\text{otherwise},
\end{cases}
$$
is odd.
\item[(b)]
$\Sel_2(E^{(\epsilon n)}/\BQ)/\kappa(E^{(\epsilon n)}(\BQ)_\tor)=0$.
\item[(c)]
$L(1,E^{(\epsilon n)})\neq 0$ and the analytic Sha
$\mathscr{L}(\epsilon n)$ of $E^{(\epsilon n)}$ is odd.
\end{enumerate}
\end{thm}

As in \cite{Tian}, the core of the approach: a link between $\mathscr{L}(\pm n)$ and $g(n)$.
In light of \cite{SD}, \cite{Kane}, \cite{Sm2}, the link also  allows us to deduce
the following positive density \cite{FLPT}.

\begin{thm}Let $E$ be the elliptic curve $y^2=x(x-1)(x+3)$.
Among the set of positive square-free integers $n\equiv 7 \pmod{24}$ (resp.~$n\equiv  3 \pmod{24}$),
the subset of $n$ - for which the analytic Sha $\mathscr{L}(\pm n)$ of $E^{(\pm n)}$ is odd and
$\Sel_2(E^{(\pm n)}/\BQ)/\kappa(E^{(\pm n)}(\BQ)_\tor)$ trivial -
has density
$\displaystyle{\frac12\prod_{i=1}^\infty(1-2^{-i})\thickapprox 14.4\%}$
(resp.~$\prod_{i=1}^\infty(1-2^{-i}) $).
\end{thm}
\begin{remark}
A salient feature: $\mathscr{L}(\pm n)$ mirrors  
the ideal class groups of imaginary quadratic fields,
unlike \cite{Kane}, \cite{Sm2}.
\end{remark}
In the following we exhibit methods to study Theorem \ref{FLPT-TFAE}.
\subsection{$2$-descent in a quadratic twist family}

Let $E/\BQ$ be an elliptic curve.

\subsubsection{The set-up}
In the following we assume $E[2]\subset E(\BQ)$.

One may suppose that $E$ is given by a Weierstrass equation
$y^2=x(x-e_1)(x-e_2)$ with $e_1,e_2\in\BZ$.
Let $m$ be a square-free integer,
and $E^{(m)}:y^2=x(x-e_1m)(x-e_2m)$.
Let $S$ be the set of primes dividing $2me_1e_2(e_1-e_2)\infty$
and $\BQ(S, 2)$ the subgroup
of $\BQ^\times/(\BQ^\times)^2$ supported on $S$.

Then elements in $\Sel_2(E^{(m)}/\BQ)$ can be  realized as curves $C_\Lambda/\BQ$ with $C_\Lambda(\BA_\BQ)\neq \emptyset$.  Here for $\Lambda=(b_1, b_2)\in \BQ(S, 2)^2$,
the curve:
$$ C_\Lambda: \qquad \begin{cases}b_1z_1^2-b_2z_2^2=e_1mt^2\\
b_1z_1^2-b_1b_2z_3^2=e_2mt^2
\end{cases}$$
Note that $C_\Lambda(\BA_\BQ)\neq \emptyset$ if and only if $C_\Lambda(\BQ_v)\neq \emptyset$ for all $v\in S$.
The $2$-torsion points $O$, $(0, 0)$, $(e_1m, 0)$, and $(e_2m, 0)$ correspond to
$(b_1, b_2)=(1, 1)$, $(e_2/e_1, -e_1m)$, $(e_1m, (e_1-e_2)/e_1)$, and $(e_2m, (e_2-e_1)m)$.
\subsubsection{The strategy}
Let $S'\subset S$ be the set of primes dividing $2e_1e_2(e_1-e_2)$,
which is independent of $m$.
Let $n>0$ be the prime-to-$S'$ part of $m$, and let $q=m/n$.
Write $n=\ell_1\cdots\ell_k$.

Let $A=(a_{ij})\in M_{k\times k}(\BF_2)$ be the R\'edei matrix associated to $n$,
defined by $a_{ij}=\left[\frac{\ell_j}{\ell_i}\right]
:=\frac12\left(1-\left(\frac{\ell_j}{\ell_i}\right)\right)$
(the additive Legendre symbol) if $i\neq j$, and $\sum_j a_{ij}=0$ for all $i$,
namely $a_{ii}=\left[\frac{n/\ell_i}{\ell_i}\right]$.
For $\Lambda=(b_1, b_2)\in \BQ(S, 2)^2$,
$t=1,2$, write $b_t=c_t\prod_i\ell_i^{x_{t,i}}$ and
$c_t=\prod_{p\in S'\cup\{-1\}}p^{y_t^{(p)}}$
for $x_t=(x_{t,1},\cdots,x_{t,k})^{\mathrm T}\in\BF_2^k$
and $y_t^{(p)}\in\BF_2$ for $p\in S'\cup\{-1\}$.

Then the condition that $C_\Lambda(\BQ_v)\neq\emptyset$ for
all $v\in S\setminus(S'\cup\{\infty\})=\{\ell_1,\cdots,\ell_k\}$
can be rephrased as a linear equation
in $\BF_2$ involving $A$ and $x_t$.
For example, consider $E:y^2=x(x-1)(x+3)$,
then the linear equation \cite{FLPT}:
$$\matrixx{A+D_q}{D_{-3}}{}{A+D_{-q}} \begin{pmatrix} x_1\\ x_2\end{pmatrix}= \begin{pmatrix} z_{c_1}\\ z_{c_2}\end{pmatrix}.$$
Here for an integer $d$ prime to $n$,
 $z_d:=\left(\left[\frac{d}{\ell_1}\right],
\cdots,\left[\frac{d}{\ell_k}\right]\right)^{\mathrm T}\in\BF_2^k$,
and $D_d:=\diag(z_d)\in M_{k\times k}(\BF_2)$.

Similarly, for each $v\in S'\cup\{\infty\}$,
the condition that $C_\Lambda(\BQ_v)\neq\emptyset$ can also be
rephrased in terms of linear algebra.
This allows us to describe $\rank_{\BF_2}\Sel_2(E^{(m)}/\BQ)$
in terms of the corank of certain ``generalized R\'edei matrix''.
Via elementary linear algebra,
one may then establish the equivalence of parts (a) and (b)
of Theorem \ref{FLPT-TFAE}.

\begin{remark}
The approach is employed in several previous works, for example,
\cite{HB2}, \cite{SD}, \cite{Kane}, \cite{Sm2}.
\end{remark}
\subsection{An induction}

Let $E/\BQ$ be an elliptic curve of conductor $N$ and
$\phi$ the associated newform.
\subsubsection{The set-up}
For a positive square-free integer $n$, let
$K=\BQ(\sqrt{-n})$ be an imaginary quadratic field,
$D$ its discriminant, and
$\eta$ the associated quadratic character.

Let $\chi:\Gal(H/K)\cong\Cl(K)\to\{\pm 1\}$ be an unramified quadratic character
over $K$, where $H$ is the Hilbert class field of $K$.
Such characters are in one-to-one correspondence with the
unramified quadratic extensions $K(\sqrt{-d})/K$,
where $d>0$ is a divisor of $n$ with $d\equiv 3\mod 4$
(resp.~$d\equiv 3\mod 4$ or $d\equiv n\mod 8$)
if $n\equiv 3\mod 4$ (resp.~$2\mid n$).
We usually denote such $\chi$ associated to $K(\sqrt{-d})/K$ by $\chi_d$.
In particular, taking $d=n$, we obtain the trivial character
$\mathbf 1_K=\chi_n$ over $K$.

Let $\Sigma$ be the set of places $v\mid N\infty$ such that
$$\varepsilon_v(E,\chi)\cdot \chi_v\eta_v(-1)=-1.$$
Here $\varepsilon_v(E,\chi)$ is the local root number at $v$ of the
Rankin--Selberg $L$-series
$L(s,E\times\chi)$.
Note that $\infty\in \Sigma$ and for any finite place $v\in \Sigma:$ $v$ non-split in $K$.
Assume that $\Sigma$ has even cardinality, equivalently, the sign of the functional equation of $L(s,E\times\chi)$ is $+1$.

Let $B$ be the definite quaternion algebra over $\BQ$ ramified exactly at the places in $\Sigma$.
Then there exists an embedding of $K$ into $B$, which we fix once and for all.
Let $R\subset B$ be an order of discriminant $N$ with $R\cap K=\CO_K$.
The Shimura set
$$X:=B^\times \bs \wh{B}^\times /\wh{R}^\times$$
is a finite set
endowed with Hecke correspondences
$T_p$ for $p\nmid N$
and $K_v^\times$-actions for $v\mid(N,D)$ by right multiplication
(cf. \cite[Lem. 3.4]{CST}).
Let $\BC[X]$ be the set of $\BC$-valued functions on $X$,
which is endowed with a natural Hermitian inner product
$\langle\ ,\ \rangle$,
Hecke operators $T_p$ for $p\nmid N$
and $K_v^\times$-actions for $v\mid(N,D)$.
Let $\BC[X]^0\subset\BC[X]$ be the orthogonal complement of
the functions on $X$ which factor through the reduced norm map
$\widehat B^\times\to\widehat\BQ^\times$.

Let $\pi$ be the cuspidal automorphic representation of $B^\times$ whose Jacquet--Langlands transfer to
$\GL_2$ is the automorphic representation generated by $\phi$.
\subsubsection{Toric periods}
The subspace
$V(\pi,\chi)$ of $\pi^{\wh{R}^\times}\subset \BC[X]^0$ where $T_p$ acts as $a_p(\phi)=a_p(E)$ for all $p\nmid N$ and $K_v^\times$ acts via $\chi_v$ for all $v\mid(N, D)$
turns out to be one-dimensional (cf. \cite{CST}). A generator of this space is referred to as a test vector.
In fact, since the coefficients of $T_p$-actions and the values of $\chi_v$ are integral,
$V(\pi,\chi)$ has a basis $f\in\BZ[X]$ of integral values.
One may choose a $\BZ$-primitive test vector $f:X\to\BZ$, i.e. the image generates $\BZ$.  Such $f$ is unique up to $\pm 1$.

The inclusion $K\ra B$ induces a map $\iota: \Cl(K)\ra X$, in view of which, define the toric period $P^0_\chi(f)$
associated to $(E, \chi)$:
$$P_\chi^0(f)=\sum_{t\in \Cl(K)} \chi(t) f(\iota(t))\in\BZ.$$
An explicit Waldspurger formula \cite{CST}. 
\begin{thm} \label{Wald} Let $E$ and $\chi$ be as above. Then
$$L(1,E\times\chi)=2^{-\mu(N, D)} \cdot \frac{8\pi^2 (\phi, \phi)_{\Gamma_0(N)}}{u^2 \sqrt{|D|}} \cdot
\frac{|P^0_\chi(f)|^2}{\pair{f, f}}.$$
Here $\mu(N, D)$ is the number of common prime factors of $N$ and $D$, $u=[\CO_K^\times: \BZ^\times]$ and $(\phi, \phi)_{\Gamma_0(N)}$
the Petersson norm of $\phi$.
\end{thm}

Now we fix a $\BZ$-primitive test vector $f\in V(\pi,\mathbf 1_K)$ for $E$
and the trivial character $\mathbf 1_K$.
Define the genus period
$$
P_0(f):=\sum_{t\in 2\Cl(K)}f(\iota(t))\in\BZ.
$$
An evident yet key fact:
if $f|_{(\widehat B^\times)^2}$ only takes odd integeral values,
then $P_0(f)\equiv g(n)\mod 2$.
Let $k=\mu(n)$ be the number of prime factors of $n$,
let $c=0$ if $n\equiv 1\mod 4$, and $c=1$ otherwise.
In particular, $\Cl(K)/2\Cl(K)\cong\Gal(H_0/K)\cong(\BZ/2\BZ)^{k-c}$
for $H_0$ the genus class field of $K$.

Then as $\chi$ varies,
\begin{equation}\label{ind1}
\sum_{\chi:\Gal(H/K)\to\{\pm 1\}}P^0_{\chi}(f)
=2^{k-c}P_0(f).
\end{equation}
\begin{remark}
For $\chi$ such that $f\notin V(\pi,\chi)$,
notice $P^0_\chi(f)=0$.
\end{remark}
\subsubsection{Induction}
In the following we assume that $E[2]\subset E(\BQ)$.
Notice $$L(s,E\times\chi)=L(s,E^{(n/d)})L(s,E^{(-d)}).$$

By considering the variation of local arithmetic invariants of elliptic curves in a quadratic twist family,
Theorem \ref{Wald} (the Waldspurger formula) 
may be rephrased as \cite{FLPT}:
\begin{equation}\label{ind2}
|P^0_\chi(f)|^2=C_{n,d}\cdot 4^{\mu(n)}u^2
\mathscr{L}(n/d)\mathscr{L}(-d)\cdot
\left(\frac{4}{\#E^{(n/d)}(\BQ)_\tor}\right)^2
\left(\frac{4}{\#E^{(-d)}(\BQ)_\tor}\right)^2.
\end{equation}
Here $\mathscr{L}(n/d)$ and $\mathscr{L}(-d)$ are the analytic Sha of $E^{(n/d)}$
and $E^{(-d)}$ in the rank $0$ case,
and $C_{n,d}$ an explicit non-zero rational constant
which depends just on the
residue classes $(n\mod M,d\mod M)$ - where $M$ an explicit constant intrinsic to $E$.

We moreover assume that $E^{(m)}(\BQ)_\tor=E^{(m)}[2]$
for all square-free integers $m\neq\pm 1$,
and that $L(1,E)\neq 0$ and
$L(1,E^{(-1)})\neq 0$.

Now, in light of \eqref{ind1} and \eqref{ind2},
for $f\notin V(\pi,\chi_1):$
$$
2^cu\left(\pm\sqrt{C_{n,n}\mathscr{L}(1)\mathscr{L}(-n)}\cdot\frac{4}{\#E(\BQ)_\tor}
+\sum_{\substack{\chi_d:\Gal(H/K)\to\{\pm 1\}\\
f\in V(\pi,\chi_d)\\
d\neq 1,n}}\pm\sqrt{C_{n,d}
\mathscr{L}(n/d)\mathscr{L}(-d)}
\right)=P_0(f).
$$
When $f\in V(\pi,\chi_1)$, an analogous formula holds.
Finally, by replacing $E$ with $E^{(-1)}$, we obtain formulae
relating $\mathscr{L}(-1)\mathscr{L}(n)$ with $\mathscr{L}(-n/d)\mathscr{L}(d)$.
Note that in these formulae, the number of prime factors of $\pm n/d$ and $\pm d$
are strictly smaller than that of $\pm n$.
This allows us to execute an induction argument on the number $k=\mu(n)$
of prime factors of $n$, proving the lower bounds for the $2$-adic valuations
and the congruence formula modulo $2$ for $\mathscr{L}(\pm n)$.

The induction leads to an equivalence of parts (a) and (c)
of Theorem \ref{FLPT-TFAE}.
Here the presence of $4$-torsion reduces the length of the recursion formula for $\mathscr{L}(\pm n)$,
which yields a very concise formula relating
$\mathscr{L}(\pm n)$ and $g(n)$, unlike \cite{TYZ}.
\begin{remark}
The induction method was introduced in \cite{Tian} and it has been employed in several previous works,
for example \cite{CLTZ}, \cite{TYZ}.
\end{remark}


\begin{thebibliography}{XXXX}
\bibitem{Ba} N. Balsam, {\em The Parity of Analytic Ranks among Quadratic Twists of Elliptic Curves over Number Fields}, thesis, Columbia 2015.

\bibitem{BKLPR} M. Bhargava, D. Kane, H. Lenstra, B. Poonen and E. Rains, {\em Modeling the distribution of ranks, Selmer groups, and
 Shafarevich--Tate groups of elliptic curves}, Cambridge J. Math. 3 (2015), 275--321.

 \bibitem{BK}S. Bloch and K. Kato, {\em $L$-functions and Tamagawa numbers of motives},
 The Grothendieck Festschrift, Vol. I, 333--400,
Progr. Math., 86, Birkh\"auser Boston, Boston, MA, 1990.
\bibitem{BSP}S. Bocherer and R. Schulze-Pillot, {\em On a theorem of Waldspurger and on Eisenstein series of Klingen type}, Math. Ann. 288 (1990), no. 3, 361--388.

\bibitem{BKO1} A. Burungale, S. Kobayashi and K. Ota, {\em Rubin's conjecture on local units in the anticyclotomic tower at inert primes}, preprint 2021.

\bibitem{BKO2} A. Burungale, S. Kobayashi and K. Ota, {\em $p$-adic $L$-functions and rational points on CM elliptic curves at inert primes}, preprint 2021.

\bibitem{BSTsv} A. Burungale, C. Skinner and Y. Tian, \emph{The Birch and Swinnerton-Dyer conjecture: a brief survey}, The Linde Hall Inaugural Math Symposium, Proc. Sympos. Pure. Math., to appear.

\bibitem{BST1} A. Burungale, C. Skinner and Y. Tian, {\em $p$-converse to a theorem of Gross--Zagier and Kolyvagin: CM elliptic curves over totally real fields}, preprint 2019.

\bibitem{BST2} A. Burungale, C. Skinner and Y. Tian, {\em Elliptic curves and Beilinson--Kato elements: rank one aspects}, preprint 2020.

    \bibitem{BST4} A. Burungale, C. Skinner and Y. Tian, {\em $p$-converse to a theorem of Gross--Zagier, Kolyvagin and Rubin, II}, in progress.


\bibitem{TA1} A. Burungale and Y. Tian, {\em $p$-converse to a theorem of Gross-Zagier, Kolyvagin and Rubin}, Invent. Math. 220 (2020), no. 1, 211--253.

\bibitem{BuTi2} A. Burungale and Y. Tian, {\em A rank zero $p$-converse to a theorem of Gross--Zagier, Kolyvagin and Rubin}, preprint 2019, submitted.


 \bibitem{CST} L. Cai, J. Shu and Y. Tian, {\em Explicit Gross--Zagier formula and Waldspurger formula},
 Algebra Number Theory 8 (2014), no. 10, 2523--2572.

\bibitem{CLTZ} J. Coates, Y. Li, Y. Tian and S.Zhai, {\em Quadratic twists of elliptic curves}, Proc. London Math. Soc. (3) 110 (2015) 357-394.

 \bibitem{CoWi} J. Coates and A. Wiles, {\em On the conjecture of Birch and Swinnerton-Dyer},
Invent. Math. 39 (1977), no. 3, 223--251.

 \bibitem{De} C. Delaunay, {\em Heuristics on Tate--Shafarevich groups of elliptic curves defined over $\BQ$}, Experiment. Math. 10(2) (2001), 191--196.
\bibitem{JD}J. Desjardins, {\em Root number of the twists of an elliptic curve}, Journal de Théorie des Nombres de Bordeaux, Tome 32 (2020) no. 1, pp. 73--101.

\bibitem{DD1}T. Dokchitser and V. Dokchitser, {\em On the Birch-Swinnerton-Dyer quotients modulo squares},
Ann. of Math. (2) 172 (2010), no. 1, 567--596.

\bibitem{DD}T. Dokchitser and V. Dokchitser,
      {\em Elliptic curves with all quadratic twists of positive rank},
      Acta Arith. 137 (2009), no. 2, 193--197.


 \bibitem{FLPT}K. Feng, Q. Liu, J. Pan and Y. Tian, {\em Toric periods and non-tiling numbers}, preprint 2021.


  \bibitem{BFH}S. Friedberg and J.Hoffstein, {\em Non-vanishing theorems for automorphic L-functions on $\GL(2)$}, Ann. of Math. (2) 142 (1995), no. 2, 385--423.

\bibitem{Go}D. Goldfeld, {\em Conjectures on elliptic curves over quadratic fields}, Number theory, Carbondale 1979. (Proc. Southern Illinois Conf., Southern Illinois Univ., Carbondale, Ill., 1979), volume 751 of Lecture Notes
 	in Math., pages 108--118. Springer, Berlin, 1979.
\bibitem{GM} D. Goldfeld and S. Munao, {\em Alexander Smith wins the first David Goss prize in number theory}, Notices Amer. Math. Soc. 66 (2019), no. 11, 1875--1878.
 \bibitem{Gr2}B. Gross, {\em Heights and the special values of L-series}. Number theory (Montreal, Que., 1985), 115–187, CMS Conf. Proc., 7, Amer. Math. Soc., Providence, RI, 1987.

     \bibitem{GZ} B. Gross and D. Zagier, {\em Heegner points and derivatives of L-series},
Invent. Math. 84 (1986), no. 2, 225--320.

 \bibitem{HTX} W. He, Y. Tian, W. Xiong, {\em Explicit theta lifting and quadratic twist L-values}, preprint.

\bibitem{HB2} D.R. Heath-Brown, {\em The size of Selmer groups for congruent number problem, II}, Inv. Math. 118, 331--370 (1994).	

 \bibitem{Hee}Heegner, K. {\em Diophantische analysis und modulfunktionenMath}, Math. Z. 56 (1952), 227--253.
	
\bibitem{JLK} J. Johnson-Leung and G. Kings, {\em On the equivariant main conjecture for imaginary quadratic fields}, J. Reine Angew. Math. 653 (2011), 75--114.

\bibitem{Kane} D. Kane, {\em On the ranks of the 2-Selmer groups of twists of a given elliptic curve}, Algebra Number Theory 7, no. 5, 1253--1279 (2013).	

\bibitem{K} K. Kato, {\em $p$-adic Hodge theory and values of zeta functions of modular forms}, Cohomologies $p$-adiques et applications arithm\'etiques. III. Ast\'erisque No. 295 (2004), ix, 117--290.

\bibitem{K1} K. Kato, {\em Iwasawa theory and generalizations}, International Congress of Mathematicians. Vol. I, 335--357, Eur. Math. Soc., Zurich, 2007.

\bibitem{KS2}N. Katz and P. Sarnak, {\em Random matrices, Frobenius eigenvalues, and monodromy},
 American Mathematical Society Colloquium Publications, 45. American Mathematical Society, Providence, RI, 1999. xii+419 pp.

\bibitem{KMR} Z. Klagsbrun, B. Mazur and K. Rubin, {\em Disparity in Selmer ranks of quadratic twists of elliptic curves},  Ann. of Math. (2) 178 (2013), no. 1, 287--320.

 \bibitem{Ko} V. Kolyvagin, {\em Euler systems}, The Grothendieck Festschrift, Vol. II, 435--483, Progr. Math., 87, Birkh\"auser Boston, Boston, MA, 1990.

\bibitem{Mon96} P. Monsky, {\em Generalizing the Birch--Stephens theorem. I. Modular curves}, Math. Z. 221 (1996), no. 3, 415--420.

     \bibitem{Mon90}P. Monsky, {\em Mock Heegner Points and Congruent Numbers}, Math. Z. 204 (1990), no. 1, 45--67.

\bibitem{PR2}B. Perrin-Riou, { \em $p$-adic $L$-functions, Iwasawa theory and Heegner points}, Bull. Soc. Math. France 115 (1987), no. 4, 399--456.

 \bibitem{P} B. Poonen, {\em Heuristics for the arithmetic of elliptic curves}, Proceedings of the International Congress of Mathematicians--Rio de Janeiro 2018. Vol. II. Invited lectures, 399--414, World Sci. Publ., Hackensack, NJ, 2018.
  \bibitem{PR}B. Poonen and E. Rains, {\em Random maximal isotropic subspaces and Selmer groups},
 J. Amer. Math. Soc. 25 (2012), no. 1, 245--269.
  	
\bibitem{Qin}H. R. Qin, {\em Congruent numbers, quadratic forms and $K_2$}, preprint 2020.

\bibitem{Ru0} K. Rubin, {\em Tate--Shafarevich groups and L-functions of elliptic curves with complex multiplication}, Invent. Math. 89 (1987), no. 3, 527--559.

\bibitem{Ru}K. Rubin, {\em The ''main conjectures'' of Iwasawa theory for imaginary quadratic fields}, Invent. Math. 103 (1991), no. 1,
25--68.	

\bibitem{Sk} C. Skinner, {\em A converse to a theorem of Gross, Zagier and Kolyvagin},  Ann. of Math. (2) 191 (2020), no. 2, 329--354.

\bibitem{SU} C. Skinner and E. Urban, {\em The Iwasawa main conjectures for $\GL_2$},   Invent. Math. 195 (2014), no. 1, 1--277.

\bibitem{Sm2}A. Smith, {\em The congruent numbers have positive natural density}, preprint,
arXiv:1603.08479.

\bibitem{Sm1}A. Smith, {\em $2^\infty$-Selmer groups, $2^\infty$-class groups, and Goldfeld's conjecture},
preprint, arXiv:1702.02325.
\bibitem{SD} P. Swinnerton-Dyer, {\em The effect of twisting on $2$-Selmer groups}, Mathematical Proceedings of the Cambridge Philosophical Society, vol. 145(2008), 513--526.
	
\bibitem{Tian} Y. Tian, {\em Congruent Numbers and Heegner Points}, Cambridge Journal of Mathematics, 2 (2014), 117-161.
\bibitem{T} Y. Tian, {\em Congruent number problem}, Proceedings of the Sixth International Congress of Chinese Mathematicians. Vol. I, 135--151, Adv. Lect. Math. (ALM), 36, Int. Press, Somerville, MA, 2017.
\bibitem{TYZ} Y. Tian, X. Yuan and S. Zhang, {\em Genus periods, genus points and congruent number problem}, Asian J. Math. 21 (2017), no. 4, 721--773.

    \bibitem{Tu}J.B. Tunnell, {\em A classical Diophantine problem and modular forms of
	weight $3/2$},  Invent. Math. 72(2), 323--334 (1983).
	
 \bibitem{Wal} J.L. Waldspurger, {\em Sur les coefficients de Fourier des formes modulaires
de poids demi-entier}, J. Math. pures et appl. 60 (1981), 375--484.

\bibitem{YZZ} X. Yuan, S.-W. Zhang and W. Zhang, {\em The Gross-Zagier formula on Shimura curves}, Annals of Mathematics Studies, vol 184. (2013) viii+272 pages.

\bibitem{Zh} W. Zhang, {\em Selmer groups and the indivisibility of Heegner points},  Cambridge Journal of Math., Vol. 2 (2014), No. 2, 191--253.
\end{thebibliography}
\end{document}